\documentclass{amsart}

\newtheorem{theorem}{Theorem}[section]
\newtheorem{lemma}[theorem]{Lemma}
\newtheorem{proposition}[theorem]{Proposition}
\newtheorem{corollary}[theorem]{Corollary}
\theoremstyle{definition}
\newtheorem{definition}[theorem]{Definition}
\newtheorem{example}[theorem]{Example}

\theoremstyle{remark}
\newtheorem{remark}[theorem]{Remark}

\numberwithin{equation}{section}

\theoremstyle{plain}

\newtheorem{problem}[theorem]{\bf Problem}

\def\int{\mathop{\roman{int}}}

\def\1{^{-1}}

\def\proof{{\bf Proof. }}
\def\endproof{\hfill \qed}

\numberwithin{equation}{section}
\begin{document}

\title[
Rips complexes and covers in the uniform category
]
   {Rips complexes and covers in the uniform category}

\author{N.~Brodskiy}
\address{University of Tennessee, Knoxville, TN 37996}
\email{brodskiy@@math.utk.edu}

\author{J.~Dydak}
\address{University of Tennessee, Knoxville, TN 37996}
\email{dydak@@math.utk.edu}
\thanks{The second-named author was partially supported
by the Center for Advanced Studies in Mathematics
at Ben Gurion University of the Negev (Beer-Sheva, Israel).}

\author{B.~Labuz}
\address{University of Tennessee, Knoxville, TN 37996}
\email{labuz@@math.utk.edu}

\author{A.~Mitra}
\address{University of Tennessee, Knoxville, TN 37996}
\email{ajmitra@@math.utk.edu}


\date{February 7, 2008.}

\begin{abstract}
James \cite{Jam} introduced uniform covering maps as an analog
of covering maps in the topological category. Subsequently
Berestovskii and Plaut \cite{BP3} introduced a theory of covers for uniform
spaces generalizing their results for topological groups \cite{BP1}-\cite{BP2}.
Their main concepts are discrete actions and pro-discrete actions, respectively.
In case of pro-discrete actions Berestovskii and Plaut provided
an analog of the universal covering space and
 their theory works well for the so-called
coverable spaces.
 As will be seen in Section
\ref{SECTION-Comparison}, \cite{BP3}
generalizes only regular covering maps in topology and pro-discrete actions may not be preserved
by compositions.
\par
In this paper we redefine the uniform covering maps and we generalize pro-discrete
actions using Rips complexes and the chain lifting property.
We expand the concept
of generalized paths of Krasinkiewicz and Minc \cite{KraMin}.
One way to do it is by embedding $X$ in a space with good local properties
and this is done in Section \ref{SECTION GenPathsInSpaces}.
Another way is by systematic
use of Rips complexes.  In the topological category one uses
paths in $X$ originating from a base point to construct the universal covering space
$\widetilde X$.
We use paths in Rips complexes
and their homotopy classes possess a natural uniform structure, a generalization
of the basic topology on $\widetilde X$.
Applying Rips complexes leads to a natural class of uniform spaces
for which our theory of covering maps works as well as the classical one, namely
the class of uniformly joinable spaces.
In the case of metric continua (compact and connected metric spaces)
that class is identical with pointed $1$-movable spaces, a well-understood
class of spaces introduced by shape theorists (see \cite{DydSeg} or \cite{MarSeg}).
 The class of pointed $1$-movable continua contains all planar
subcontinua (examples: Hawaiian Earring and the suspension of the Cantor set) and is preserved
by continuous maps.
The most notable continuum not being pointed $1$-movable is the dyadic solenoid.
As an application of our results we present an exposition in \cite{BDLM}
of Prajs' \cite{Pra} homogeneous curve that is path-connected but not locally connected.
\end{abstract}
\maketitle

\medskip
\medskip
\tableofcontents

\section{Introduction}

The aim of this paper is to develop a theory of covering maps
in the uniform category via generalizations of the classical
construction of universal covering spaces.
For basic facts on uniform spaces we refer to \cite{Isb} or \cite{Jam}.

In Section \ref{SECTION covers} we provide an analog of covering maps in topology
adopted for the uniform category.
Our definition uses local structure of the base space just as it does in topology.
However, we provide a characterization of uniform covering maps via chain lifting
property and that characterization is later on expanded to define generalized
uniform covering maps.

How to construct universal covering space for uniform spaces $X$ with good local properties
(the so-called uniform Poincare spaces)?
Let us recall briefly the construction of a simple topology (used in \cite[p.82]{Spa},
\cite[p.253]{HilWyl}, \cite{FisZas}, \cite{BDLMPeano}, and \cite{BDLM2})
on the space $\widetilde X$, the space of homotopy classes (rel.endpoints)
of paths in $X$ originating from the base-point $x_0$.
First, one defines sets $B([\alpha],U)$
(denoted by $<\alpha,U>$ on p.82 in \cite{Spa}), where $U$ is open in $X$, $\alpha$ joins $x_0$
and $\alpha(1)\in U$ as follows: $[\beta]\in B([\alpha],U)$ if and only if there is a path $\gamma$
in $U$ from $\alpha(1)$ to $\beta(1)$ such that  $\beta\sim \alpha\ast\gamma$.

$\widetilde X$ equipped with the topology (which we call the
{\bf basic topology on $\widetilde X$})
whose basis consists of $B([\alpha],U)$, where $U$ is open in $X$, $\alpha$ joins $x_0$
and $\alpha(1)\in U$ is denoted by $\widehat X$
as in \cite{BogSie}.

It turns out, for uniform spaces $X$, the space $\widetilde X$ has a natural uniform structure that generalizes the basic topology and we provide natural analogs of classical
results for uniform Poincare spaces.

\par How to deal with spaces $X$ whose local structure is complicated
(example: the Topologist Sine Curve)? Spaces like that may not be path-connected
resulting in the projection $\widetilde X\to X$ not being surjective. The geometrical answer is
to use paths in neighborhoods of $X$. That leads to the concept
of a {\bf generalized path} introduced by Krasinkiewicz-Minc \cite{KraMin}.
We generalize that concept to embeddings of $X$ in a space $T$
with good local properties in
Section \ref{SECTION GenPathsInSpaces}. The resulting space $GP_T(X,x_0)$
of generalized paths has a natural uniform structure mimicking that of $\widetilde X$.
The advantage of embeddings is that
many natural spaces are defined that way and we may apply shape-theoretical results. The disadvantage of defining
universal covering spaces using only embeddings is that one has to show
independence of the construction on the embedding. That is why Rips complexes
are useful. In Section \ref{SECTION gen unif
paths} we apply Rips complexes to define an abstract space $GP(X,x_0)$ of
generalized uniform paths equipped with a natural uniform structure so that the
end-point map $\pi_X\colon GP(X,x_0)\to X$ is uniformly
continuous. As the defining characteristic of covering maps we use
 the Unique Path Lifts
Property of any topological cover:
 $f\colon X\to Y$ is declared {\bf a generalized uniform covering map} (see Section \ref{SECTION gen unif covers}) if it has lifting and approximate uniqueness of lifts properties for both chains
 and generalized uniform paths.
The meaning of our definition is that not only we want the Unique Generalized Path Lifting
Property but the lifting function ought to be a morphism in the uniform category.
\par
What is the largest class of spaces for which that definition ought to work?
The answer is quite simple:
it is the class of {\bf uniformly joinable spaces} $X$ that may be characterized
by the requirement
of $\pi_X\colon GP(X,x_0)\to X$ being a generalized uniform covering map.
It turns out that particular class (in case of metric continua) coincides with the class
of joinable continua studied by Krasinkiewicz and Minc \cite{KraMin}.

In Section \ref{SECTION-Comparison} we relate our construction to that of Berestovskii and Plaut \cite{BP3}.

We are grateful to Conrad Plaut for a series of lectures on his work with Berestovskii.
We thank Misha Levin for suggesting to provide an exposition of J.Prajs' \cite{Pra}
example of a homogeneous curve $P$  that is path-connected but not locally connected
(see \cite{BDLM}).

\section{Uniform covering maps}\label{SECTION covers}

We will discuss exclusively symmetric subsets $E$ of $X\times X$
(that means $(x,y)\in E$ implies $(y,x)\in E$) and the natural notation here (see \cite{Pla})
is to use $f(E)$ for the set of pairs $(f(x),f(y))$, where $f\colon X\to Y$ is a function.
Similarly, $f^{-1}(E)$ is the set of pairs $(x,y)$ so that $(f(x),f(y))\in E$ if $f\colon X\to Y$
and $E\subset Y\times Y$.
\par The {\bf ball $B(x,E)$ at $x$ of radius $E$} is the set of all $y\in X$
satisfying $(x,y)\in E$.
\par
A {\bf uniform structure} on $X$ is a family $\mathcal{E}$ of symmetric subsets $E$
of $X\times X$ (called {\bf entourages}) that contain the diagonal
of $X\times X$, form a filter (that means $E_1\cap E_2\in \mathcal{E}$
if $E_1,E_2\in \mathcal{E}$ and $F_1\in \mathcal{E}$ if $F_2\in \mathcal{E}$
and $F_2\subset F_1$),
and every $G_1\in \mathcal{E}$ admits $G\in \mathcal{E}$ so that
$G^2\subset G_1$ ($G^2$ consists of pairs $(x,z)\in X\times X$
so that there is $y\in X$ satisfying $(x,y)\in G$ and $(y,z)\in G$).
A {\bf base} $\mathcal{F}$ of a uniform structure $\mathcal{E}$
is a subfamily $\mathcal{F}$ of $\mathcal{E}$ so that for every entourage $E$ there is a subset
$F\in \mathcal{F}$ of $E$.

\par Given a decomposition of a uniform space $X$ the most pressing issue is
if it induces a natural uniform structure on the decomposition space. James \cite[2.13 on p.24]{Jam}
has a concept of weakly compatible relation to address that issue. For the purpose of this paper
we need a different approach.
\begin{definition}\label{GeneratingUCSTructureDef}
Suppose $f\colon X\to Y$ is a surjective function from a uniform space
$X$. {\bf $f$ generates a uniform structure on $Y$} if the family
$f(E)$, $E$ an entourage of $X$, is a base of a uniform structure on $Y$
(that particular uniform structure on $Y$ is said to be {\bf generated by $f$}).
Equivalently, for each entourage $E$ of $X$ there is an entourage $F$ of $X$
such that $f(F)^2\subset f(E)$.
\end{definition}

Notice $f\colon X\to Y$ is uniformly continuous if both $X$ and $Y$ are uniform spaces
and the uniform structure on $Y$ is generated by $f$. Indeed $E\subset f^{-1}(f(E))$
for any entourage $E$ of $X$.

Uniform covering maps were defined by James \cite[p.112]{Jam}.
In this section we redefine that concept using Rips complexes and we provide a characterization
of uniform covering maps in terms of chain lifting. That characterization will be very useful when
generalizing uniform covering maps in Section \ref{SECTION gen unif covers}.

The definition of a Rips complex for uniform structures is a straightforward
generalization of Rips complexes \cite[Chapter 4]{GhyHar} for metric spaces.

\begin{definition}\label{DefRipsComplex}
Let $X$ be a set.
Given a symmetric subset $E$ of $X\times X$ containing the diagonal define
the {\bf Rips complex} $R(X,E)$ as the subcomplex of the full complex over $X$
whose simplices are finite subsets $F=\{x_1,\ldots,x_n\}$ of $X$ so that $F\times F\subset E$.
\end{definition}

Notice $E$ containing the diagonal of $X\times X$ ensures the set of vertices of $R(X,E)$ coincides with $X$.

Given $f\colon X\to Y$ and an entourage $E$ of $X$
notice it induces a natural simplicial map $f_E\colon R(X,E)\to R(Y,f(E))$
by the formula $f_E(\sum\limits_{i=1}^n t_i\cdot x_i)=\sum\limits_{i=1}^n t_i\cdot f(x_i)$.

Our goal is to study homotopy classes of paths in $R(X,E)$ joining two of
its vertices. Since the identity function $K_w\to K_m$, $K$ a
simplicial complex, from $K$ equipped with the CW (weak) topology
to $K$ equipped with the metric topology is a homotopy equivalence
(see \cite[page 302]{MarSeg}), it does not really matter which
topology we choose for $R(X,E)$. For simplicity (and to be able to use
\cite[Corollary 17 on p.138]{Spa}), let it be the
weak topology.

The simplest path in $R(X,E)$ is the {\bf edge-path} $e(x,y)$ starting from
$x$ and ending at $y$ so that $(x,y)\in E$.

Any path in $R(X,E)$ joining two vertices $x$ and $y$ can be
realized, up to homotopy (see \cite[Section 3.4]{Spa}), as a
concatenation of edge-paths. Thus, each path in $R(X,E)$ can be
realized by an {\bf $E$-chain} $x_1=x$, \ldots, $x_n=y$ such that
$(x_i,x_{i+1})\in E$ for all $i< n$. Two paths in $R(X,E)$
represented by different $E$-chains with the same end-points are
homotopic rel. end-points if and only if one can move from one
chain to the other by simplicial homotopies: a new vertex $v$ can
be added or removed from a chain if and only if $v$ forms a
simplex in $R(X,E)$ with adjacent links of a chain (see
\cite[Section 3.6]{Spa}).

Here is our definition of covering maps in the uniform category using
Rips complexes. We call a simplicial map a {\bf simplicial covering map}
if it is a topological cover.

\begin{definition}\label{DefOfECover}
Let $X$ and $Y$ be uniform spaces.
 $f\colon X\to Y$ is a
{\bf uniform covering map} if it generates the uniform structure on $Y$ and
the family $\mathcal{E}$ of entourages  $E$
of $X$ such that the induced map $f_E\colon R(X,E)\to R(Y,f(E))$
is a simplicial covering map forms a base of the uniform structure of $X$.
\end{definition}

Let us characterize uniform covering maps in terms
analogous to classical topological covering maps.

\begin{definition}\label{DefOfEvenlyCovered}
Let $f\colon X\to Y$ be a map of sets. A symmetric subset $E$ of
$X\times X$ {\bf evenly covers} $f(E)$
 if $B(x,E)$ is mapped by $f$ bijectively onto $B(f(x),f(E))$ for all $x\in X$.
\end{definition}

\begin{lemma}\label{JamesVsRipsDefinitionOfUCM}
Suppose $f\colon X\to Y$ is a function of uniform spaces
and the uniform structure on $Y$ is generated by $f$.
$f\colon X\to Y$
is a uniform covering map if and only if $X$ has a base of entourages
$E$ that evenly cover $f(E)$.
\end{lemma}
\proof Suppose $f_E\colon R(X,E)\to R(Y,f(E))$
is a simplicial covering map. If $(x,y), (x,z)\in E$ and $f(y)=f(z)$, then the edge-path $e(f(x),f(y))$
can be lifted starting from $x$ in two different ways unless $y=z$.
That means $B(x,E)$ is mapped by $f$ injectively into $B(f(x),f(E))$ for all $x\in X$.
If $(f(x),y)\in f(E)$ we can lift the edge $e(f(x),y)$ to an edge $e(x,z)$ in $R(X,E)$.
Thus $f(z)=y$ and $(x,z)\in E$.
\par Suppose $B(x,E)$ is mapped by $f$ bijectively onto $B(f(x),f(E))$ for all $x\in X$.
Assume $E^2\subset F$ and $F$ covers evenly $f(F)$.
 Given $x\in X$ and given a simplex $\Delta$ in $R(Y,f(E))$
containing $f(x)$ there is a unique lift $\Delta^{\prime}$ of $\Delta$
containing $x$. Indeed, we can lift each edge of $\Delta$ emanating from $f(x)$
and the endpoints of lifts (together with $x$) form a lift $\Delta^{\prime}$.
If $[f(x),y_1]$ and $[f(x),y_2]$ are two edges of $\Delta$
that lift to $[x,x_1]$ and $[x,x_2]$ respectively then there is an edge
$[x_1,z]$ of $R(X,E)$ with $f(z)=y_2$. Now $x_2,z\in B(x_1,F)$ so $x_2=z$.
Thus every open star of a vertex in $R(Y,f(E))$ has the point inverse
of the form of the disjoint union of open stars of vertices in $R(X,E)$
and $f_E$ restricts to a homeomorphism on each of those open stars.
In other words, $f_E$ is a topological cover.
\endproof

To show how \ref{DefOfECover} relates to uniform covering maps
of James \cite[p.112]{Jam} let us define one of the main concepts of
the paper.
\begin{definition}\label{ChainLiftingProperty}
A surjective function $f\colon X\to Y$ from a uniform space $X$
has {\bf the chain lifting property} if for any entourage $E$ of $X$ there is an entourage
$F$ of $X$ such that any $f(F)$-chain in $Y$ starting from $f(x_0)$
can be lifted to an $E$-chain starting from $x_0$.
\par A function $f\colon X\to Y$ from a uniform space $X$
has {\bf the uniqueness of chain lifts property} if for every entourage $E$ of $X$
there is an entourage $F$ of $X$ such that any two $F$-chains
$\alpha$ and $\beta$ satisfying $f(\alpha)=f(\beta)$ are equal if they
originate from the same point.
\end{definition}

Notice the chain lifting property is stronger than generating a uniform structure
on the range (see \ref{ChainLiftingPropertyAndUniformOpeness}).

James \cite[p.13]{Jam} defined the concept of an entourage
being transverse to an equivalence relation. In the same way one can
define an entourage to be transverse to a function.

\begin{definition}\label{TransverseEntourage}
Let $X$ be a uniform space and $Y$ be a set.
An entourage $E$ of $X$ is {\bf transverse} to $f\colon X\to Y$ if
$(x,y)\in E
$ and $f(x)=f(y)$ implies $x=y$.
\end{definition}

\begin{proposition}\label{UniquenessOfChainLiftsAndTransverseEnt}
$f\colon X\to Y$ has the uniqueness of chain lifts property if and only if
$f$ has a transverse entourage.
\end{proposition}
\proof Suppose $E$ is an entourage of $X$ and $F\subset E$
is chosen so that $F^2$ is transverse to $f$.
Given two different $F$-chains $\alpha=\{x_0,\ldots,x_n\}$ and $\beta=\{y_0,\ldots,y_n\}$ of $X$ originating from $x_0$
such that $f(\alpha)=f(\beta)$ choose the smallest $i$ satisfying $x_i\ne y_i$.
Notice $(x_i,y_i)\in F^2$ as $x_{i-1}=y_{i-1}$. Hence $x_i=y_i$ (as $F^2$ is transverse
to $f$), a contradiction.
\par If $f$ has the uniqueness of chain lifts property, pick an entourage $E_0\subset X\times X$
so that any two $E_0$-chains $\alpha$ and $\beta$ are equal provided $f(\alpha)=f(\beta)$
and their origins are the same. If $f(x)=f(y)$ and $(x,y)\in E_0$,
then put $\alpha=\{x,y\}$ and $\beta=\{x,x\}$. Observe $f(\alpha)=f(\beta)$.
Hence $\alpha=\beta$ and $E_0$ is transverse to $f$.
\endproof

Here is the relation of chain lifting property to the concept
of uniform openness used by James \cite[Definition 1.12, p.10]{Jam}.

\begin{proposition}\label{ChainLiftingPropertyAndUniformOpeness}
Suppose $f\colon X\to Y$ is a surjective function from a uniform space
$X$ to a set $Y$.
\begin{itemize}
\item[a.] If $Y$ has a uniform structure making $f$ uniformly open,
then $f$ has the chain lifting property.
\item[b.] If $f$ has the chain lifting property, then
$f$ generates a uniform structure on $Y$
making $f$ uniformly open.
\end{itemize}
\end{proposition}
\proof a. $f$ being uniformly open means existence, for each entourage
$D$ of $X$, of an entourage $E$ of $Y$ such that $B(f(x),E)\subset f(B(x,D))$.
That condition says any pair $(y,z)\in E$ lifts to $(x,t)\in E$ if $y=f(x)$.
Hence any $E$-chain in $Y$ lifts to a $D$-chain in $X$.
Choose an entourage $F$ of $X$ satisfying $f(F)\subset
E$
and notice any $f(F)$-chain lifts to
a $D$-chain.

b. Suppose $f$ has the chain lifting property.
First, we need to show the family $\{f(E)\}_{E\in\mathcal{E}(X)}$
forms a base of entourages of $Y$.
The only condition needed to be proved is the existence,
for each entourage $E$ of $X$, of an entourage $F$ of $X$ such that $f(F)^2\subset
f(E)$.
Assume $D^2\subset E$
and any $f(F)$-chain in $Y$ lifts to a $D$-chain in $X$.
Suppose $(f(x),y)\in f(F)$ and $(y,z)\in f(F)$.
We may choose $x_1\in f^{-1}(y)$ so that $(x,x_1)\in D$. Now we may choose $x_2\in f^{-1}(z)$
so that $(x_1,x_2)\in D$. Hence $(x,x_2)\in E$ and $(f(x),z)\in f(E)$.
\par Notice that if any $f(F)$-chain in $Y$ lifts to an $E$-chain in $X$,
then $B(f(x),f(F))\subset f(B(x,E))$, so $f$ is indeed uniformly open.
\endproof

Let us characterize covering maps in the uniform category in terms of lifting of chains.

\begin{theorem}\label{CoversInTermsOfChainLifting}
Suppose $f\colon X\to Y$ is a function of uniform spaces
and the uniform structure on $Y$ is generated by $f$.
$f$ is a uniform covering map if and only if
the following conditions are satisfied:
\begin{enumerate}
\item[a)] $f$ has the chain lifting property.
\item[b)] $f$ has the uniqueness of chain lifts property.
\end{enumerate}
\end{theorem}
\proof Suppose $f$ is a uniform covering map.
The existence of a transverse entourage $E_0$ is obvious as any $E$ such that $f_E\colon R(X,E)\to R(Y,f(E))$
is a simplicial covering map will do.
Condition b) follows from~\ref{UniquenessOfChainLiftsAndTransverseEnt}.
Also, in that case it is clear any $f(E)$-chain in $Y$ lifts
to an $E$-chain in $X$.
\par
Assume Conditions a) and b) are satisfied.
\par
Given an entourage $G$ of $X$ define $\alpha(G)$
as the set of points $(x,y)\in G$ satisfying the following conditions:
\begin{enumerate}
\item For any $x_1\in f^{-1}(f(x))$ there is $y_1\in f^{-1}(f(y))$
such that $(x_1,y_1)\in G$.
\item For any $y_1\in f^{-1}(f(y))$ there is $x_1\in f^{-1}(f(x))$
such that $(x_1,y_1)\in G$.
\end{enumerate}
First, observe the family $\{\alpha(G)\}_{G\in\mathcal{E}}$ forms a base of entourages of $X$.
Indeed, given an entourage $E$ choose an entourage $E_1\subset E$
so that any $f(E_1)$-chain in $Y$ lifts to an $E$-chain in $X$. Now $E_1\subset \alpha(E)$
as follows: given $(x,y)\in E_1$ and given $x_1\in f^{-1}(f(x))$,
one can lift the $f(E_1)$-chain $f(x),f(y)$ to an $E$-chain $x_1,y_1$. Similarly,
if $y_1\in f^{-1}(f(y))$, we can lift $f(y),f(x)$ to an $E$-chain $y_1,x_1$.
That means $(x,y)\in \alpha(E)$.
\par
Second,
if $E_0$ is an entourage of $X$ transverse to $f$ (provided by~\ref{UniquenessOfChainLiftsAndTransverseEnt}), then
notice $\alpha(G)^2\subset E_0$ implies $f$ maps $B(x,\alpha(G))$
bijectively onto $B(f(x),f(\alpha(G)))$. Indeed, if $y,z\in B(x,\alpha(G))$
and $f(y)=f(z)$, then $(y,z)\in E_0$ and $y=z$.
If $z\in B(f(x),f(\alpha(G)))$, there is $(x_1,y_1)\in \alpha(G)$ so that $f(x)=f(x_1)$
and $f(y_1)=z$. As $(x_1,y_1)\in \alpha(G)$ there must exist $y\in f^{-1}(f(y_1))$
satisfying $(x,y)\in G$. Notice that implies $(x,y)\in \alpha(G)$ (as $(x_1,y_1)\in \alpha(G)$)
and that means any $f(\alpha(G))$-chain lifts to an $\alpha(G)$-chain.
\endproof

\begin{remark}\label{ComparisonOfOurDefAndJames}
James \cite[p.111--112]{Jam} defined uniform covering maps
as $p\colon X\to B$
so that there is an entourage $E$ transverse to $p$
and $X$ has a base of entourages $F$ satisfying
$R\circ F=F\circ R$, where $R=p^{-1}(\Delta B)$ is the relation on $X$
induced by $p$. Unfortunately, he never added the condition
that $p$ generates the uniform structure on $Y$ (in the language of \cite{Jam} it translates to
relation $R$ being weakly compatible with the uniform structure on $X$). Our interpretation of Chapter 8 in \cite{Jam} is that
he assumes so implicitly. With that in mind our definition
of uniform covering maps is equivalent to that of James.
Indeed, \ref{UniquenessOfChainLiftsAndTransverseEnt}
takes care of the uniqueness of chain lifts property and \ref{ChainLiftingPropertyAndUniformOpeness} implies that James' uniform covering
maps have the chain lifting property as a map that generates the uniform structure and statisfies
$R\circ F=F\circ R$ is uniformly open. Conversely, observe that any $F$ that evenly covers $p(F)$
satisfies $R\circ F=F\circ R$.
\end{remark}

\par
The most important property of covering maps in topology is that of unique lifts of
paths and the fact homotopic paths have the same end-point when lifted.
That leads to a quick candidate $\widetilde X$ for the universal cover of a pointed
space $(X,x_0)$: it is the quotient space of the space of paths $Map((I,0),(X,x_0))$ in $X$
(equipped with the compact-open topology)
starting from $x_0$, where the equivalence relation is that of homotopy
rel. end-points.

In the reminder of this section we are going to define a uniform structure
on $\widetilde X$ mimicking the basic topology on $\widetilde X$ and we are going to discuss
necessary and sufficient conditions for the projection
$\pi_X\colon \widetilde X\to X$ ($\pi_X(\alpha)$ is the end-point of $\alpha$)
to be a uniform covering map. It turns out, not surprisingly,
those conditions involve uniform local path-connectedness
and uniform semi-local simple connectedness. However, our definition
of uniform local path-connectedness is much simpler than
\cite[Definition 8.12 on p.119]{Jam} and we are unsure if the definition
\cite[Definition 8.13 on p.119]{Jam} of uniform semi-local simple connectedness
is correct as it involves existence of a base of entourages rather than just one entourage.

For each entourage $E$ of $X$
define $E^\ast$ as the family of pairs of homotopy classes
$([\alpha],[\beta])$ of paths from $x_0$ such that
$\alpha^{-1}\ast \beta$ is homotopic rel. end-points to a path
contained in some $B(z,E)$. The family $\{E^\ast\}_E$ forms a
base of a uniform structure on $\widetilde X$
which we call the {\bf basic uniform structure} on $\widetilde X$. Notice that the
projection $\pi_X\colon \widetilde X\to X$
is uniformly continuous. Also, it is surjective if and only if $X$ is path-connected.

\begin{proposition}\label{UnifLocalPathConnectCharacterization}
If $X$ is a path-connected uniform space, then
its structure is generated by
$\pi_X\colon
\widetilde X\to X$ if and only if
for each entourage $E$ of $X$ there is an entourage $F$
such that any two points in $B(x,F)$ can be connected by a path in
$B(x,E)$ for any $x\in X$.
\end{proposition}
\proof Suppose $\pi_X\colon
\widetilde X\to X$ is generates the structure on $X$. Given an entourage $E$ of $X$
there is an entourage $F$ of $X$ satisfying $F^2\subset \pi_X(E^\ast)$.
Suppose $y,z\in B(x,F)$. Since $(y,z)\in F^2$, there is a pair of paths
$(\alpha,\beta)\in E^\ast$ so that $\alpha$ joins $x_0$ to $y$ and
$\beta$ joins $x_0$ to $z$. Thus $\alpha^{-1}\ast \beta$ is homotopic
rel. end-points to a path contained in some $B(w,E)\subset B(x,E^3)$.
\par
Suppose for each entourage $E$ of $X$ there is an entourage $F$
such that any two points in $B(x,F)$ can be connected by a path in
$B(x,E)$ for any $x\in X$. Given $(y,z)\in F$
choose a path $\alpha$ contained in $B(y,E)$ joining $y$ to $z$.
Choose a path $\beta$ from $x_0$ to $y$
and observe $(\beta,\beta\ast\alpha)\in E^\ast$,
$\pi_X(\beta,\beta\ast\alpha)=(y,z)$ which proves
$F\subset \pi_X(E^\ast)$ and $\pi_X$ generates the structure on $X$.
\endproof

Call a space satisfying the conditions in \ref{UnifLocalPathConnectCharacterization}
{\bf uniformly locally path-connected} (see \cite{Pla} or \cite{BP3}).

\begin{lemma}
\label{UnifLocPathConnLemma}
If $X$ is uniformly locally path-connected, then for every entourage $E$ of $X$ there is
an entourage $F\subset E$ such that $B(x,F)$ is path-connected for every $x\in X$.
\end{lemma}

\proof Let $H$ be an entourage of $X$ such that for any $x\in X$, any two points in $B(x,H)$
can be joined by a path in $B(x,E)$. Define $F$ to be all $(x,y)\in E$ such that $x$ and $y$
can be joined by a path in some $B(z,E)$. Notice $H\subset F$ so that $F$ is an entourage.
Let $x\in X$ and $y,z\in B(x,F)$. Then there is a path $\alpha$ joining $y$ to $x$ in some
$B(z_1,E)$ and a path $\beta$ joining $x$ to $z$ in some $B(z_2,E)$. Notice that $\alpha \ast \beta$
is contained in $B(x,F)$. \endproof

\begin{proposition} \label{LiftsCovInUnifCat}
Suppose $f:X\rightarrow Y$ is a uniform covering map and
$g\colon Z\rightarrow Y$ is uniformly continuous. Suppose $X,Y$, and $Z$ are
path-connected and uniformly locally path-connected. Let $x_{0}\in X$,
$y_{0}\in Y$, and $z_{0}\in Z$ with $f(x_{0})=g(z_{0})=y_{0}$. Then there is
a unique uniformly continuous lift $\widetilde{g}:Z\rightarrow X$ of $g$ with
$\widetilde{g}(z_{0})=x_{0}$ if and only if $g_{\ast }(\pi _{1}(Z,z_{0}))\subset
f_{\ast }(\pi _{1}(X,x_{0}))$. Further, if $g$ is a uniform
covering map then so is $\widetilde{g}$.
\end{proposition}

\proof From the corresponding theorem in the topological category we have
the forward direction. Also from the topological theorem there is a unique
lift $\widetilde{g}$ with $\widetilde{g}(z_{0})=x_{0}$ defined by letting $\widetilde{g}
(z)$ be the endpoint of the lift of $g\circ \alpha $ starting at $x_{0}$
where $\alpha $ is a path from $z_{0}$ to $z$ \cite{Mun}. We show that
$\widetilde{g}$ is uniformly continuous. Let $E$ be an entourage of $X$
evenly covering $f(E)$ and $F$ be an entourage of $Z$ such that $B(z,F)$ is
path-connected for each $z\in Z$ and $g(F)\subset f(E)$. Let $(x,y)\in F$.
Take a path $\alpha $ from $z_{0}$ to $x$ and a path $\beta $ from $x$ to $y$
that is contained in $B(x,F)$. Lift $g\circ \alpha $ to a path $\widetilde{\alpha
}$ from $x_{0}$ to $\widetilde{g}(x)$. Since $g\circ \beta $ is contained in $
B(g(x),f(E))=B(f(\widetilde{g}(x)),f(E))$, $\gamma =(f|_{B(\widetilde{g}
(x),E)})^{-1}\circ g\circ \beta $ is a path starting at $\widetilde{g}(x)$ that is
contained in $B(\widetilde{g}(x),E)$. But $\widetilde{\alpha}\ast \gamma $ is a lift
of $g\circ (\alpha \ast \beta )$ starting at $x_{0}$ so $\gamma (1)=\widetilde{g}
(y)$ and we have $(\widetilde{g}(x),\widetilde{g}(y))\subset E$.

Now suppose $g$ is a uniform covering map. Let us first
show that $\widetilde{g}$ is surjective. If $x\in X$, take a path $\alpha $ from
$x_{0}$ to $x$ and lift $f\circ \alpha $ to a path $\widetilde{\alpha}$ in $Z$
starting at $z_{0}$. Then $\widetilde{g}(\widetilde{\alpha}(1))=x$.

Now let us see that $\widetilde{g}$ generates the uniform structure of $X$. Suppose $E$ is an
entourage of $Z$ evenly covering $g(E)$ and let $F$ be an entourage of $X$
that evenly covers $f(F)$ and has $f(F)\subset g(E)$. Finally, let $G\subset
F$ be an entourage so that for every $x\in X$, any two points in $B(x,G)$
can be joined by a path contained in $B(x,F)$. Suppose $(x,y)\in G$. Take a
path $\alpha $ from $x_{0}$ to $x$ and a path $\beta $ from $x$ to $y$ that
is contained in $B(x,F)$. Lift $f\circ \alpha $ to a path $\widetilde{\alpha}$
in $Z$ starting at $z_{0}$ and set $x^{\prime }=\widetilde{\alpha}(1)$. Set $
\gamma =(g|_{B(x^{\prime },E)})^{-1}\circ f\circ \beta $ and notice that $\widetilde{
\alpha}\ast \gamma $ is a path from $z_{0}$ to some $y^{\prime }\in
B(x^{\prime },E)$. Then $\widetilde{g}(x^{\prime },y^{\prime })=(x,y)$ so we
have $G\subset \widetilde{g}(E)$ and $\widetilde{g}$ the uniform structure of $X$.

Finally, put $H=E\cap \widetilde{g}^{-1}(G)$ where $E$ and $G$ are as above and
let us see that $H$ evenly covers $\widetilde{g}(H)$. Let $z\in Z$ and suppose $
x,y\in B(z,H)$ with $\widetilde{g}(x)=\widetilde{g}(y)$. Then $g(x)=g(y)$ so $x=y$
since $E$ evenly covers $g(E)$. Now suppose $y\in B(\widetilde{g}(z),\widetilde{g}
(H))$. Take a path $\alpha $ in $Z$ from $z_{0}$ to $z$ and let $\widetilde{
\alpha}$ be the lift of $g\circ \alpha $. Now take a path $\beta $ from $
\widetilde{g}(z)$ to $y$ that is contained in $B(\widetilde{g}(z),F)$. Set $\gamma
=(g|_{B(z,E)})^{-1}\circ f\circ \beta $ and notice $\alpha \ast \gamma $ is a
path from $z_{0}$ to some $y^{\prime }\in B(z,E)$. Then $\widetilde{g}(y^{\prime
})=y$ and $y^{\prime }\in B(z,H)$.\endproof

When is $\pi_X\colon \widetilde X\to X$ a uniform covering map?

\begin{proposition}\label{Semi-LocallySimplyConnectedChar}
Let $X$ be a path-connected uniform space.
Suppose $E$ is an entourage of $X$ and $x\in X$.
If $(E^\ast)^2$ is transverse to $\pi_X$, then every loop in $B(x,E)$ at $x$ is null-homotopic
in $X$.
\end{proposition}
\proof Notice balls $B(\alpha,E^\ast)$, $\alpha\in \pi_X^{-1}(x)$, are mutually
disjoint. Suppose $\gamma$ is a loop in $B(x,E)$ at $x$.
Choose $\alpha$ joining $x_0$ to $x$ and notice
$(\alpha\ast\gamma,\alpha)\in E^\ast$. Since $\pi_X(\alpha)=x=\pi_X(\alpha\ast\gamma)$,
$\alpha\ast\gamma$ is homotopic rel. end-points to $\alpha$ in $X$
and $\gamma$ is null-homotopic in $X$.
\endproof

\ref{UnifLocalPathConnectCharacterization}
and \ref{Semi-LocallySimplyConnectedChar}
 lead to the concept of a {\bf uniform Poincare
space} $X$ (compare \cite{BP3}), a space that is path-connected,
uniformly locally path-connected, and {\bf uniformly semi-locally
simply connected} (that means the existence of an entourage $F$
such that all loops in $B(x,F)$ at $x$ are null-homotopic in $X$ for all
$x\in X$).

\begin{theorem}\label{PoincareSpacesAndUnifCovers}
$\pi_X\colon \widetilde X\to X$
is a uniform covering map if and only if $X$ is a uniform
Poincare space.
\end{theorem}
\proof If $\pi_X$ is a uniform covering map,
$X$ must be uniformly locally path-connected by \ref{UnifLocalPathConnectCharacterization}
and uniformly semi-locally simply connected by \ref{Semi-LocallySimplyConnectedChar}.

Suppose $X$ is a uniform Poincare space. By \ref{UnifLocalPathConnectCharacterization}
$\pi_X$ generates the uniform structure of $X$. Let $F$ be an entourage of $X$ and let $E$ be an entourage of $X$ such that
loops in $B(x,E)$ at $x$ are null-homotopic in $X$. Let $G$ be an entourage
with $G^{2}\subset F\cap E$ and $H\subset G$ be an entourage such that all
balls $B(x,H)$ are path-connected
(use~\ref{UnifLocPathConnLemma}).

Let us show that $H^{\ast }$ evenly covers $\pi _{X}(H^{\ast })$. Let
$\alpha \in \widetilde{X}$ and $\beta ,\gamma \in B(\alpha ,H^{\ast })$ with
$\pi_{X}(\beta )=\pi _{X}(\gamma )=x$ for some $x\in X$. Notice
$\beta ^{-1}\ast\gamma $ is
homotopic rel. end-points to a path
contained in $B(x,H^{2})\subset B(x,E)$ so it is null homotopic.
Therefore $\beta \sim \gamma$. Now let $y\in B(\pi _{X}(\alpha ),\pi_{X}(H^{\ast }))$,
so $y=\pi _{X}(\beta )$ and $\pi _{X}(\alpha )=\pi_{X}(\alpha ^{\prime })$
for some $(\beta ,\alpha ^{\prime })\in H^{\ast }$.
Then $y,\pi _{X}(\alpha )\in B(z,H)$ for some $z\in X$ so there is a path $\gamma$
joining them that is contained in $B(z,H)$. Notice $\pi _{X}(\alpha \ast \gamma )=y$
and $(\alpha ,\alpha \ast \gamma )\in H^{\ast }$. \endproof

\section{Generalized uniform paths}\label{SECTION gen unif paths}

How to adjust the above construction of $\widetilde X$
for spaces with bad local properties?
A good way is to approximate $X$ by its Rips complexes.
An alternative way is to embed $X$ in a space with good local
properties and use paths there (see Section \ref{SECTION GenPathsInSpaces}).

First, we will extend the concept of paths being homotopic rel. end-points.

\begin{definition}\label{DefOfEHomotopy}
Two paths $c$ and $d$ in $R(X,E)$
with endpoints in $X$
are {\bf $E$-homotopic}
provided the following conditions are satisfied:

\begin{itemize}
\item[1.] The initial points $x_c$ and $x_d$ and the terminal points $y_c$ and $y_d$ of the
paths $c$ and $d$ satisfy $(x_c,x_d),(y_c,y_d)\in E$;
\item[2.] $c$ is homotopic in $R(X,E)$ rel. end-points to the concatenation
$e(x_c,x_d)\ast d\ast e(y_d,y_c)$.
\end{itemize}
\end{definition}

Notice the relation of being $E$-homotopic is symmetric and
coincides with usual homotopy of paths rel. end-points in $R(X,E)$ if the
end-points of paths are the same.

Given a uniform space $X$ one can consider the space $GP(X)$ of
generalized paths in $X$. A {\bf generalized path} is a collection
$\{[c_E]\}_E$ of homotopy classes of paths $[c_E]$ in $R(X,E)$
joining fixed $x\in X$ to $y\in X$ such that for all entourages
$F\subset E$, $c_F$ is homotopic to $c_E$ in $R(X,E)$ rel.
end-points.
\par
A generalized path $c=\{[c_E]\}_E$ is called {\bf $F$-short} if
its end-points $x$ and $y$ satisfy $(x,y)\in F$ and $[c_F]$ is the
homotopy class of the edge-path $e(x,y)$ in $R(X,F)$. In other
words, $c$ is $F$-short if $c_F$ is $F$-homotopic to the constant
path at the origin of $c$.
\par
We equip $GP(X)$ with a natural uniform structure: a base of entourages of $GP(X)$ is the family $F^\ast$
consisting of all pairs $(c,d)$ of generalized
paths $c=\{[c_E]\}_E$ and $d=\{[d_E]\}_E$
such that $c_F$ is $F$-homotopic to $d_F$.

If two generalized paths $c$ and $d$ have the same initial point (or the same terminal point),
then $(c,d)\in F^\ast$ if and only if $c^{-1}\ast d$ is $F$-short ($c\ast d^{-1}$ is $F$-short,
respectively).

\par The {\bf projection} $\pi_X\colon GP(X)\to X$
assigns to each generalized path its end-point. Notice $\pi_X$
is uniformly continuous as $(c,d)\in F^\ast$ implies $(\pi_X(c),\pi_X(d))\in F$.
\par

Given a uniform morphism $f\colon X\to Y$
it induces a function $\widetilde f\colon GP(X)\to GP(Y)$ as follows:
Let $c=\{[c_E]\}_E$ be any generalized path of $X$ and $F$ be any entourage of $Y$.
Put $E=f^{-1}(F)$ and define $\widetilde f(c)=\{[f_E(c_E)]\}_F$. Notice that all paths $f_H(c_H)$, $H\subset E$,
are homotopic rel. end-points to $f_E(c_E)$ so that $\widetilde f$ is well-defined.
Also notice that $\widetilde f$ is uniformly continuous as for any entourages $E$ of $X$ and $F$
of $Y$ the inclusion $E\subset f^{-1}(F)$ implies $\widetilde f(E^\ast)\subset F^\ast$.
\par
Given a pointed uniform space $(X,x_0)$ one can consider
the space $GP(X,x_0)$ of generalized paths in $X$ originating from $x_0$ with the uniform structure induced from $GP(X)$. Any pointed uniformly continuous
function $f\colon (X,x_0)\to (Y,y_0)$ induces a uniformly continuous
$\widetilde f\colon GP(X,x_0)\to GP(Y,y_0)$.

In case of $(X,x_0)=(I,0)$ being the pointed unit interval the space $GP(I,0)$
is naturally identical with $I$ as for any $t\in I$ there is only one
generalized uniform path from $0$ to $t$
(a generalization of this observation is Corollary~\ref{GenPathsForUnifPoincare}).
Therefore every ordinary path in $X$ from $x_0$ to $x$
induces naturally a generalized uniform path which we will usually denote
by the same letter.

\section{Uniform joinability}\label{Uniform joinability}

Connectivity and path connectivity can be generalized to the uniform category
in several ways.
First, the concept of {\bf chain connectivity} of $X$ (see \cite{Pla} or \cite{BP3})
that is equivalent to {\bf uniform connectivity} of James \cite[Definition 1.5 on p.7]{Jam}
 can be formulated
as connectivity of all its Rips complexes.
\par Here is a generalization of path-connectivity.

\begin{definition}\label{DEfJoinability}
$X$ is {\bf joinable} if any of its two points
can be joined by a generalized uniform path.
\end{definition}

Obviously, any $X$ such that the underlying topological space
is path-connected, is joinable.

The following is an elementary exercise:
\begin{proposition}\label{JoinabilityEqSurjectivity}
If $X$ is a uniform space, then the following conditions are equivalent:
\begin{itemize}
\item[a.] $X$ is joinable,
\item[b.] $\pi_X\colon GP(X,x_0)\to X$ is surjective for each $x_0\in X$,
\item[c.] $\pi_X\colon GP(X,x_0)\to X$ is surjective for some $x_0\in X$.
\end{itemize}
\end{proposition}

\begin{definition}\label{DEfUJ}
$X$ is {\bf uniformly joinable} if for each entourage $E$ of $X$
there is an entourage $F$ such that any pair $(x,y)\in F$
can be joined by a generalized path $\{[c_H]\}_H\in GP(X)$
that is $E$-short.
\end{definition}

Notice that any uniformly locally path-connected
$X$
 is uniformly joinable. Those include inner-metric spaces (in particular, geodesic spaces)
 and Peano continua.

\begin{proposition}\label{Bi_UImageOfUJIsUJ}
If $f\colon X\to Y$ generates the uniform structure of $Y$ and $X$ is uniformly joinable, then
$Y$ is uniformly joinable.
\end{proposition}
\proof Given an entourage $E$ of $Y$ pick an entourage $F\subset
G=f^{-1}(E)$ of $X$ so that for any pair $(x,y)\in F$ there is a
generalized path $c(x,y)$ joining $x$ and $y$ that is $G$-short.
Suppose $(x^{\prime},y^{\prime})\in f(F)$. Pick a pair $(x,y)\in
F$ satisfying $f(x,y)=(x^{\prime},y^{\prime})$ and observe $\widetilde
f(c(x,y))$ is a generalized path in $Y$ joining $x^{\prime}$ and
$y^{\prime}$ whose $E$-th term is $e(x^{\prime},y^{\prime})$ in
$R(Y,E)$.
\endproof

\begin{proposition}\label{ChainConnectedAndUJAreJoinable}
If $X$ is uniformly joinable and chain connected, then
it is joinable.
\end{proposition}
\proof Given an entourage $E$ of $X$ pick an entourage $F$
of $X$ so that any pair $(y,z)\in F$ can be connected by a generalized
path $c(y,z)$ that is $E$-short.
Since $x_0$ and $x_1$ can be connected by an $F$-chain, we can replace
each link of that chain by a generalized path and obtain a generalized path $d$
from $x_0$ to $x_1$.
\endproof

\begin{proposition}\label{BiUniformAndJoinability}
If $X$ is chain connected, then the following conditions are equivalent:
\begin{itemize}
\item[a.] $X$ is uniformly joinable,
\item[b.] $\pi_X\colon GP(X,x_0)\to X$ generates the uniform structure of $X$ for each $x_0\in X$,
\item[c.] $\pi_X\colon GP(X,x_0)\to X$ generates the uniform structure of $X$ for some $x_0\in X$.
\end{itemize}
\end{proposition}
\proof a)$\implies$b).
$\pi_X$ is surjective
by~\ref{ChainConnectedAndUJAreJoinable}
and~\ref{JoinabilityEqSurjectivity}. Given an entourage $E$ of $X$
pick an entourage $F$ of $X$ so that any pair $(y,z)\in F$ can be
connected by a generalized path $c(y,z)$ that is $E$-short. Let
$d$ be a generalized uniform path from $x_0$ to $y$. Now, $d\ast
c(y,z)$ is a generalized path in $X$ so that $(d,d\ast c(y,z))\in
E^\ast$. Since $\pi_X(d)=y$ and $\pi_X(d\ast c(y,z))=z$, we obtain
$F\subset \pi_X(E^\ast)$ which proves $\pi_X$ generates the uniform structure of $X$.

c)$\implies$a).
If $\pi_X$ generates the uniform structure of $X$, then for each entourage $E$ of $X$ there is
an entourage $F\subset E$ of $X$ such that $F\subset \pi_X(E^\ast)$. That means
for any pair $(x,y)\in F$ there is $(c,d)\in E^\ast$ with $x=\pi_X(c)$ and $y=\pi_X(d)$.
Notice $e=c^{-1}d$ is a generalized $E$-short path from $x$ to $y$. \endproof

\begin{definition}\label{DefOfFundamentalUProGroup}
Suppose $X$ is a uniform space and $x_0\in X$.
By the {\bf uniform fundamental pro-group} $\mathrm{pro-}\pi_1(X,x_0)$
we mean the inverse system of groups $\{\pi_1(R(X,E),x_0)\}_E$.
\par The {\bf uniform fundamental group} $\check\pi_1(X,x_0)$
is the inverse limit of $\mathrm{pro-}\pi_1(X,x_0)$ which is identical
with the group of generalized loops of $X$ at $x_0$.
Notice $\check\pi_1(X,x_0)$ inherits a uniform structure from $GP(X,x_0)$,
so it is actually a topological group.
\end{definition}

Recall an inverse system $\{G_a\}_{a\in A}$ of groups
satisfies the {\bf Mittag-Leffler condition} (see \cite[p.77]{DydSeg} or \cite[p.165]{MarSeg}) if for every $a\in A$ there is $b > a$ such that for any $c > b$ the image of $G_b\to G_a$ is contained in the image of $G_c\to G_a$
(that implies those images are actually equal).
In particular, an inverse system $\{G_a\}_{a\in A}$ of groups is {\bf trivial}
if for every $a\in A$ there is $b > a$ such that the image of $G_b\to G_a$ is trivial.

As noted in \cite[Proposition 6.1.2]{DydSeg} an inverse system of groups
satisfies the Mittag-Leffler condition if and only if it is movable
in the category of pro-sets. Therefore it makes sense to consider
a condition equivalent to uniform movability (see \cite[p.160]{MarSeg}) of a pro-group
in the category of pro-sets.

\begin{definition}\label{StrongMLDefinition}
An inverse system $\{G_a\}_{a\in A}$ of groups with inverse limit $G$
satisfies the {\bf strong Mittag-Leffler condition} if for every $a\in A$ there is $b > a$ such that
the image of $G\to G_a$ contains the image of $G_b\to G_a$
(that implies those images are actually equal).
\end{definition}

\begin{remark}\label{countable ML implies strong ML}
If the index set $A$ has a countable cofinal set $B$ (that means for any $a\in A$
there is $b\in B$ with $b\ge a$), then $\{G_a\}_{a\in A}$ satisfying Mittag-Leffler condition
implies it satisfying the strong Mittag-Leffler condition (use \cite[Theorem 4 on p.163]{MarSeg}).
\end{remark}

\begin{theorem}\label{UJImpliesStrongML}
If $X$ is uniformly joinable
then $\mathrm{pro-}\pi_1(X,x)$ satisfies the strong Mittag-Leffler condition for each $x\in X$.
\end{theorem}
\proof Fix $x\in X$. Given an entourage $E$ of $X$ pick an
entourage $F\subset E$ with the property that any pair of points
$(y,z)\in F$ can be connected by a generalized path $c(y,z)$ so
that $c(y,z)_E$ is the homotopy class of the edge $e(y,z)$ in
$R(X,E)$. Suppose $\alpha$ is a loop at $x$ in $R(X,F)$. Represent
that loop as an $F$-chain $x=x_1,\ldots, x_n=x$ and replace each
edge $e(x_i,x_{i+1})$ by $c(x_i,x_{i+1})$. The result is a
generalized loop $\gamma$ at $x$ so that $[\gamma_E]=[\alpha]$ in
$R(X,E)$. Notice that $\gamma_H$, $H\subset F$, represents an
element of $\pi_1(R(X,H),x)$ whose image in $\pi_1(R(X,E),x)$ is
the same as $[\alpha]$.
\endproof

\begin{theorem}\label{StrongMLPlusJoinImpliesUJ}
Suppose $X$ is a joinable uniform space.
If $\mathrm{pro-}\pi_1(X,x_0)$ satisfies the strong Mittag-Leffler condition for some $x_0\in X$,
then $X$ is uniformly joinable.
\end{theorem}
\proof Given an entourage $E$ of $X$ choose an entourage $F\subset E$
with the property $im(\pi_1(R(X,F),x_0)\to \pi_1(R(X,E),x_0))\subset im(\check\pi_1(X,x_0)\to \pi_1(R(X,E),x_0))$. If $(x,y)\in F$ choose a generalized path $c(x)$ from $x_0$ to $x$
and choose a generalized path $c(y)$ from $x_0$ to $y$.
The loop $c(x)_F\ast e(x,y)\ast c(y)_F^{-1}$ in $R(X,F)$ equals $d_E$ in $R(X,E)$
for some generalized loop $d$ at $x_0$.
Consider $c=c(x)^{-1}\ast d\ast c(y)$ and notice $c_E=e(x,y)$
 which proves $X$
is uniformly joinable.
\endproof

 \begin{proposition}\label{Lim1AndJoinability}
If $X$ is a chain connected uniform space, then any of the
following conditions implies that $X$ is joinable:
\begin{itemize}
\item[(1)] $\mathrm{pro-}\pi_1(X,x_0)$ is trivial for some $x_0\in
X$; \item[(2)] $X$ has a countable base of entourages and
${\lim\limits_{\leftarrow}}^1 (\mathrm{pro-}\pi_1(X,x_0))=0$ for
some $x_0\in X$. \end{itemize}
\end{proposition}
\proof (1) Consider the Rips complex $R(\mathcal{E})$ of the family of entourages of $X$.
A simplex in $R(\mathcal{E})$ is a finite set $\Delta=\{E_1,\ldots,E_k\}$
of entourages of $X$ such that for any pair $i,j\leq k$ either $E_i\subset E_j$
or $E_j\subset E_i$. Each $\Delta$ has a minimal vertex $m(\Delta)$ defined
as $\bigcap\limits_{i=1}^kE_i$. By induction on the number of vertices of
$\Delta$ find an entourage $a(\Delta)\subset m(\Delta)$
such that the inclusion-induced
homomorphism $\pi_1(R(X,a(\Delta)),x_0)\to\pi_1(R(X,m(\Delta)),x_0)$ is
trivial and the function $a(\Delta)$ is monotone
(i.e. $a(\Delta)\subset a(\Delta^{\prime})$ whenever $\Delta^{\prime}$
is a face of $\Delta$).
\par
Fix a point $x\in X$.
Then any two
paths from $x_0$ to $x$ in $R(X,a(E))$ are homotopic in $R(X,E)$.
Let $c_E$ be such a path. Then $\{[c_E]\}_E$ is a generalized path
from $x_0$ to $x$. Indeed, if $F\subset E$, then for $\Delta=\{F,E\}$
one has $a(\Delta)\subset a(F)\subset F$ and $a(\Delta)\subset a(E)\subset E$,
so a path in $R(X,a(\Delta))$ from $x_0$ to $x$
is homotopic rel. end-points to both $c_E$ and $c_F$ in $R(X,E)$.

(2) Let $E_n$ be a base of entourages of the uniform structure on
$X$. We can assume $E_{i+1}\subset E_i$ for all $i$ since $\bigcap\limits_{i=1}^nE_i$
is also a base. Put $G_n=\pi_1(R(X,E_n),x_0)$. Recall (see \cite{MarSeg})
that ${\lim\limits_{\leftarrow}}^1\{G_n\}=0$ means
existence, for each sequence $g_n\in G_n$, of a sequence $h_n\in
G_n$ such that $g_k=h_k\cdot h_{k+1}^{-1}$ in $G_k$ for all $k\ge
1$. If each $G_n$ is countable, that condition is equivalent to
$\{G_n\}$ satisfying the Mittag-Leffler condition (see \cite[p.78]{DydSeg}).

Given $x\in X$ choose, for each $n\ge 1$, a path $p_n$ in
$R(X,E_n)$ from $x_0$ to $x$. Put $g_n=[p_n\ast p^{-1}_{n+1}]$ and
choose loops $h_n$ at $x_0$ so that $h_k\ast h_{k+1}^{-1}$ is
homotopic rel. $x_0$ to $g_k$ in $R(X,E_k)$. Put
$c_k=p_k^{-1}\ast h_k$ for $k\ge 1$. For each $E$ choose $E_k\subset E$
and set $c_E=c_k$. We then have a generalized path in $X$ from $x_0$ to $x$.
\endproof

\begin{proposition}\label{PathUniqueness}
If $X$ is a uniformly joinable uniform space, then
\par\noindent $\mathrm{pro-}\pi_1(GP(X,x_0),y_0)$ is trivial for any $x_0\in X$,
where $y_0$ is the constant generalized path at $x_0$ in $X$.
\end{proposition}

\proof Given an entourage $E$ of $X$ choose an entourage $F\subset E$ such
that any two points $(x,y)\in F$ can be connected by a generalized path $
c(x,y)$ that is $E$-short. Take any loop in $R(GP(X),F^{\ast })$ based at $
y_{0}$ and represent it as a sequence $y_{0},\ldots ,y_{k}=y_{0}$ of
generalized paths in $X$. Let $x_{i}$ be the endpoint of $y_{i}$. Then $
x_{0},\ldots ,x_{k}$ is an $F$-chain that is $F$-homotopic to $(y_{1}\ast
y_{1}^{-1}\ast y_{2}\ast y_{2}^{-1}\ast \cdots \ast y_{k-1}\ast
y_{k-1}^{-1})_{F}$ and is therefore null-homotopic via a finite sequence of
simplicial homotopies in $R(X,F)$. We wish to mimick those simplicial
homotopies in $R(GP(X),E^{\ast }).$ At each stage of the homotopy we will
have an $E^{\ast }$-chain in $GP(X,x_{0})$ such that the endpoints of the
links of the chain form an $F$-chain in $X.$ In case of a vertex reduction,
say $x_{i}$, the sequence $y_{0},\ldots ,y_{i-1},y_{i+1},\ldots ,y_{k}$ is
an $E^{\ast }$-chain since $(y_{i-1}^{-1}\ast y_{i+1})_{E}$ is homotopic to $
(y_{i-1}^{-1}\ast y_{i}\ast y_{i}^{-1}\ast y_{i+1})_{E}$ which in turn is $E$
-homotopic to $e(x_{i-1},x_{i})\ast e(x_{i},x_{i+1})$ and the simplex $
[x_{i-1},x_{i},x_{i+1}]\in R(X,E)$. Also the endpoints $x_{0},\ldots
,x_{i-1},x_{i+1},\ldots x_{k}$ form an $F$-chain. In the case of inserting a
new vertex $z$ between $x_{i}$ and $x_{i+1}$ we create a new sequence $
y_{0},\ldots ,y_{i},y_{i}\ast c(x_{i},z),y_{i+1},\ldots ,y_{k}$. This
sequence is an $E^{\ast }$-chain since $((y_{i}\ast c(x_{i},z))^{-1}\ast
y_{i+1})_{E}$ is $E$-homotopic to $e(z,x_{i})\ast e(x_{i},x_{i+1})$ and the
simplex $[z,x_{i},x_{i+1}]\in R(X,E)$. Again, the endpoints form an $F$-chain.
\endproof

\begin{corollary}\label{GPIsUJ}
If $X$ is uniformly joinable, then $GP(X,x_0)$ is chain connected
and uniformly joinable for any $x_0\in X$.
\end{corollary}
\proof Put $Y=GP(X,x_0)$ and let $y_0$ be the constant generalized
path at $x_0$ in $X$. Given an entourage $E$ of $X$ choose an
entourage $F\subset E$ such that any pair $(x,y)\in F$ can be
connected by a generalized path $c(x,y)$ whose $E$-th term is
$e(x,y)$. If $c$ is an element of $GP(X,x_0)$ look at $c_{F}$ and
pick its simplicial representative, an edge-path $x_0$,
$x_1$,\ldots, $x_n$. Let $d$ be the concatenation of
$c(x_i,x_{i+1})$, $i=0,\ldots, n-1$. Put $e=c\ast d^{-1}$ and
notice $(y_0,e)\in E^\ast$. Now the sequence $y_0$, $y_1=e$,
$y_2=e\ast c(x_0,x_1)$, \ldots, $y_{n+1}=c$ (here $y_{i+1}=y_i\ast
c(x_{i-1},x_i)$) joins $y_0$ and $c$ so that $(y_i,y_{i+1})\in
E^\ast$ for all $i$. Thus $Y$ is chain connected. Application
of~\ref{PathUniqueness},~\ref{Lim1AndJoinability},
and~\ref{StrongMLPlusJoinImpliesUJ}
completes the proof.
\endproof

\begin{corollary}\label{GPSquare=GP}
If $X$ is uniformly joinable, then for any $x_0\in X$ the
projection $\pi_{GP(X,x_0)}\colon GP(GP(X,x_0),c)\to GP(X,x_0)$ is
a uniform equivalence for any $c\in GP(X,x_0)$.
\end{corollary}
\proof By \ref{GPIsUJ} the space $Y=GP(X,x_0)$ is chain connected
and uniformly joinable. By \ref{BiUniformAndJoinability}
$\pi_Y\colon GP(Y,c)\to Y$ generates the uniform structure of $Y$ and by
\ref{PathUniqueness} it is injective. Therefore $\pi_Y$ is a
uniform equivalence.
\endproof

\section{Generalized uniform covering maps}\label{SECTION gen unif covers}

We define generalized uniform covering maps
by weakening conditions of \ref{CoversInTermsOfChainLifting}
(for relations between uniform covering maps and generalized uniform
covering maps via inverse limits see \cite{Lab}).

\begin{definition}\label{DefGenUnifCov}
 A {\bf generalized uniform covering map} is a function
$f\colon X\to Y$ of uniform spaces generating the uniform structure of $Y$
and satisfying the following conditions:
\begin{enumerate}
\item[GP1.] ({\bf Generalized Path Lifting Property}) Every generalized uniform path in $Y$ at $f(x_0)$ lifts to
a generalized uniform path in $X$ at $x_0$.
\item[GP2.] ({\bf Approximate Uniqueness of Generalized Path Lifts Property})
For any entourage $E$ of $X$ there is an entourage $F$ of $X$ such that any two generalized
uniform paths $\alpha$ and $\beta$ in $X$ with a common origin must be $E$-homotopic if $f(\alpha)$ and $f(\beta)$ are $f(F)$-homotopic.

\item[C1.] $f$ has the chain lifting property.
\item[C2.]
For any entourage $E$ of $X$ there is an entourage $F$ of $X$ such that any two $F$-chains $\alpha$ and $\beta$ with a common origin must be $E$-homotopic if $f(\alpha)$ and $f(\beta)$ are $f(F)$-homotopic.
\end{enumerate}
\end{definition}

Notice that Conditions C1 and C2 are discrete versions of Conditions GP1 and GP2, respectively.

Before analyzing interdependence of Conditions GP1-2 and C1-2 let us
explain the meaning of Conditions GP1-2.

\begin{proposition}\label{GUCMInduceUHomeo}
Suppose $X$ and $Y$ are Hausdorff uniform spaces and $f\colon X\to Y$ generates the uniform structure of $Y$.
\begin{enumerate}
\item If $f$ satisfies Conditions GP1-2, then $\widetilde f\colon GP(X,x_0)\to GP(Y,f(x_0))$
is a uniform equivalence for each $x_0\in X$.
\item If $X$ is joinable and $\widetilde f\colon GP(X,x_1)\to GP(Y,f(x_1))$
is a uniform equivalence for some $x_1\in X$,
then $f$ satisfies Conditions GP1-2.
\end{enumerate}
\end{proposition}
\proof (1) Condition GP1 of \ref{DefGenUnifCov} says $\widetilde f$ is surjective
and Condition GP2 of \ref{DefGenUnifCov} implies $\widetilde f$ is both injective and generates the uniform structure of $GP(Y,f(x_0))$. Indeed, if $\widetilde f(\alpha)=\widetilde f(\beta)$, then $\alpha$ is $E$-homotopic
to $\beta$ for all entourages $E$ of $X$. Hence their end-points coincide and $\alpha=\beta$.
Condition GP2 means (provided GP1 holds) $F^\ast\subset \widetilde f(E^\ast)$,
so $\widetilde f$ generates the uniform structure of $GP(Y,f(x_0))$.
\par (2) Suppose $\alpha$ is a generalized uniform path in $Y$ starting
at $f(x_0)$. Choose a generalized uniform path $\gamma$ from $x_1$ to $x_0$
and let $\beta$ be a generalized uniform path from $x_1$ satisfying $\widetilde f(\beta)=\widetilde f(\gamma)\ast \alpha$. Put $\sigma=\gamma^{-1}\ast \beta$ and observe $\widetilde f(\sigma)=\alpha$.
That proves GP1.
\par
Choose an entourage $F$ of $Y$ so that $F^\ast \subset \widetilde f(E^\ast)$
(such $F$ exists as $\widetilde f$ is a uniform equivalence).
Suppose $\alpha$ and $\beta$ are two generalized uniform paths
at $x_0$ such that $f(\alpha)$ is $F$-homotopic to $f(\beta)$.
Choose
a generalized uniform path $\gamma$ from $x_1$ to $x_0$
and observe $(\widetilde f(\gamma\ast \alpha),\widetilde f(\gamma\ast\beta))\in F^\ast$.
That implies there are two generalized uniform paths
$(\alpha_1,\beta_1)\in E^\ast$ starting from $x_1$ so that $\widetilde f(\alpha_1)=\widetilde f(\gamma\ast \alpha)$ and $\widetilde f(\beta_1)=\widetilde f(\gamma\ast \beta)$. Due to $\widetilde f$ being injective,
$\alpha_1=\gamma\ast \alpha$ and $\beta_1=\gamma\ast \beta$. Now $\alpha^{-1}\ast\beta=\alpha_1^{-1}\ast\beta_1$
is $E$-short and Condition GP2 holds.
\endproof

\begin{corollary}\label{UniquenessOfLifts}
Suppose $f\colon X\to Y$ is a generalized uniform covering map
and $y_0=f(x_0)$.
If $Z$ is joinable, then for every uniformly continuous $g\colon Z\to Y$
with $g(z_0)=y_0$
there is at most one uniformly continuous lift $h\colon Z\to X$ of $g$ satisfying $h(z_0)=x_0$.
\end{corollary}
\proof Given $z\in Z$ pick a generalized path $c$ from $z_0$ to $z$.
Since $\widetilde f(\widetilde h(c))=\widetilde g(c)$, the generalized path $\widetilde h(c)$ is uniquely
determined. Hence its end-point $h(z)$ is uniquely determined as well.
\endproof

The following result has an easy proof, so it is left to the reader.

\begin{proposition}\label{CompositionOfGenCovers}
Suppose $f\colon X_1\to X_2$ and $g\colon X_2\to X_3$ generate the uniform structure of their ranges.
If $f$ and $g$ are generalized uniform covering maps, then so is the composition $g\circ f$.
\end{proposition}

Our next objective is to replace Condition C2 by Approximate Uniqueness of Chain Lifts Property
as it is closer to the uniqueness of lifts property in our definition of uniform covering maps.

\begin{lemma}
\label{TwoConditionsOnGUCMImplyTheThird}
 Given a function $f\colon X\to Y$ from a  uniform space $X$ consider the following conditions:

\begin{itemize}
\item[C2.] For any entourage $E$ of $X$ there is an entourage $F$ of $X$ such that any two $F$-chains $\alpha$ and $\beta$ with a common origin must be $E$-homotopic if $f(\alpha)$ and $f(\beta)$ are $f(F)$-homotopic.
\item[C3.] ({\bf Approximate Uniqueness of Chain Lifts Property})
For any entourage $E$ of $X$ there is an entourage $F$ of $X$ such that any
two $F$-chains in $X$ starting from $x_0$ are $E$-close if their images are identical.
\end{itemize}
If $f$ has the chain lifting property, then C2 and C3 are equivalent.
\end{lemma}
\proof C2$\implies$C3. If $f(\alpha)=f(\beta)$, then they are clearly $f(F)$-homotopic.
Hence $\alpha$ is $E$-homotopic to $\beta$. In particular, their end-points are $E$-close.
The same argument works of subchains of $\alpha$ and $\beta$ with the same number of links,
so $\alpha$ is $E$-close to $\beta$.
\par C3$\implies$C2. Given an entourage $E$ of $X$ choose an entourage $E_a$
satisfying $E_a^2\subset E$. Now choose an entourage $E_b$ of $X$ so that
any two $E_b^2$-chains must be $E_a$-close if their images are identical.
 Pick an entourage $F\subset E_b$ of $X$ such that any
$f(F)$-chain in $Y$ lifts to an $E_b$-chain in $X$.
\par Consider two $F$-chains $\alpha$ and $\beta$ starting from $x_0$ with common end-point
such that
$f(\alpha)$ and $f(\beta)$ are $f(F)$-homotopic rel.end-points. Let $\gamma_1$,\ldots,$\gamma_n$ be a sequence of $f(F)$-chains
realizing $f(F)$-homotopy from $f(\alpha)$ to $f(\beta)$. Choose an $E_b$-lift $\lambda_i$ of $\gamma_i$
for each $1< i < n$ and put $\lambda_1=\alpha$, $\lambda_n=\beta$.
To show $\lambda_i$ is $E$-homotopic to $\lambda_{i+1}$ it suffices to consider the case
$\gamma_{i+1}$ is obtained from $\gamma_i$ via an $f(F)$-expansion.
Let $\Gamma$ be the chain obtained from $\lambda_{i+1}$ by dropping the expansion vertex.
Notice $\Gamma$ is an $E_b^2$-chain and $f(\Gamma)=f(\lambda_i)$.
Hence $\Gamma$ is $E_a$-close to $\lambda_i$ and it is $E_a^2$-homotopic rel.end-points to
$\lambda_i$. Since $\lambda_{i+1}$ is an $E_b$-expansion of $\Gamma$,
it is $E$-homotopic to $\lambda_i$.
\par We still need to show that the endpoints of $\alpha$ and $\beta$ are sufficiently close. Pick an entourage $F_1\subset F$ of $X$ so that any $f(F_1)$-chain in $Y$ lifts to an $F$-chain in $X$.
Assume $\alpha$ and $\beta$ are $F_1$-chains starting from $x_0$ such that
$f(\alpha)$ and $f(\beta)$ are $f(F_1)$-homotopic.
Extend $f(\alpha)$ to $\mu$ by adding the end-point of $f(\beta)$ and lift $\mu$ to
an $F$-chain $\alpha^\prime$. Now $f(\alpha^\prime)$ and $f(\beta)$ are $f(F)$-homotopic,
so by the previous case $\alpha^\prime$ is $E$-homotopic to $\beta$ rel.end-points.
Since $\alpha^\prime$ with end-point removed is $E_a$-close to $\alpha$,
we get $\alpha$ is $E^2$-homotopic to $\beta$.
\endproof

\ref{TwoConditionsOnGUCMImplyTheThird}
says the difference between \ref{CoversInTermsOfChainLifting}  and Conditions C1-2 of
\ref{DefGenUnifCov}
is that for uniform covering maps one has existence and uniqueness of lifts
of chains (assuming the chains are sufficiently fine - that comes from existence
of an entourage transverse to the covering map)
and for generalized uniform covering maps one has existence
and approximate uniqueness of lifts of chains.

Let us show that Condition GP2 is superfluous  and
 Condition GP1 in \ref{DefGenUnifCov}
follows from C1 and C2 provided the fibers of $f$ are complete.

\begin{proposition}\label{GP1andC1-2InduceGP2}
If $f\colon X\to Y$ satisfies Conditions GP1 and C1-2 of  \ref{DefGenUnifCov}, then it satisfies Condition GP2 of \ref{DefGenUnifCov}.
\end{proposition}
\proof
Suppose any two $F$-chains $\alpha$ and $\beta$ originating from the same point
 must be $E$-homotopic if $f(\alpha)$
and $f(\beta)$ are $F_1$-homotopic, where $F_1=f(F)$. Given $(\widetilde f(d_1),\widetilde f(d_2))\in F_1^\ast$,
where $d_1,d_2\in GP(X,x_0)$, we may assume $F$-terms of both $d_1$ and $d_2$
are realized by $F$-chains $\alpha$ and $\beta$, respectively.
Since $f(\alpha)$ is $F_1$-homotopic to $f(\beta)$,
$\alpha$ is $E$-homotopic to $\beta$ resulting in $(d_1,d_2)\in E^\ast$.
\endproof

\begin{lemma}\label{CompleteFibersImplyConditionC}
Suppose $f\colon X\to Y$ is a uniformly continuous map
with complete fibers.
$f$ is a generalized uniform covering map if it generates the uniform structure on $Y$
and
the Conditions C1-2 of \ref{DefGenUnifCov} are satisfied.
\end{lemma}
\proof
The only task is to show that any generalized path $c=\{[c_{E}]\}$ in $Y$
starting at $f(x_{0})$ has a lift starting at $x_{0}$. Given an entourage $E$
of $X,$ choose an entourage $\alpha (E)\subset E$ so that $f\alpha (E)$-chains in $Y$ lift to $E$-chains. Let $\widetilde{c_{f\alpha (E)}}$ be an $E $-lift of $c_{f\alpha (E)}$ and define $x_{E}$ to be the endpoint of $\widetilde{c_{f\alpha (E)}}.$ To see that $\{x_{E}\}$ is Cauchy let $E$ be
an entourage of $X$ and choose an entourage $H$ of $X$ so that two $H$-chains are $E$-homotopic if their images are $f(H)$-homotopic. Suppose $F_{1},F_{2}\subset H.$ Then $\widetilde{c_{f\alpha (F_{1})}}$ and $\widetilde{c_{f\alpha (F_{2})}}$ are $H$-chains with images $c_{f\alpha
(F_{1})}$ and $c_{f\alpha (F_{2})}$ respectively. Since $f\alpha
(F_{1}),f\alpha (F_{2})\subset f(H),$ $c_{f\alpha (F_{1})}$ and $c_{f\alpha
(F_{2})}$ are $f(H)$-homotopic so $(x_{F_{1}},x_{F_{2}})\in E.$ Let $x$ be a
limit point of $\{x_{E}\}.$ We will extend lifts of chains $c_{E}$ to $x$ to
form a generalized path $d=\{[d_{E}]\}$ so $f(d)=c.$
\par Given an entourage $E$ of $X$ choose an entourage $\beta (E)\subset E$ so
that two $\beta (E)$-chains are $E$-homotopic if their images are $f\beta (E)
$-homotopic and choose an entourage $\gamma (E)\subset \beta (E)$ so that
for any entourage $F\subset \gamma (E),$ $(x_{F},x)\in \beta (E)$. Given an
entourage $E$ of $X$ define $d_{E}$ to be $\widetilde{c_{f\alpha \gamma (E)}}
$ extended to $x.$ Since $\widetilde{c_{f\alpha \gamma (E)}}$ and $d_{E}$
are both $\beta (E)$-chains with $f\beta (E)$-homotopic images, $f(d_{E})=c_{f\alpha \gamma (E)}$ which is $E$-homotopic to $c_{f(E)}$ so $f(d)=c.$ To see that $d$ is in fact a generalized path suppose $F\subset E$
are entourages of $X$ and consider the entourage $H=f\alpha (\alpha \gamma
(F)\cap \alpha \gamma (E)).$ Choose an $\alpha \gamma (F)\cap \alpha \gamma
(E)$-lift $h$ of $c_{H}$ and notice it is a $\beta (F)$-chain whose image is
$f\beta (F)$-homotopic to $c_{f\alpha \gamma (F)}.$ Therefore $h$ is $F$-homotopic to $d_{F}.$ Similarly $h$ is $E$-homotopic to $d_{E}$ so we have $d_{F}$ $E$-homotopic to $d_{E}.$
\endproof

\subsection{Generalized uniform covering maps and uniformly joinable spaces}

\begin{theorem}\label{AltCharOfGenUnifCovForUJ}
Suppose $X$ is uniformly joinable chain connected Hausdorff uniform space.
If $f\colon X\to Y$ generates the uniform structure of $Y$, then the following conditions are
equivalent:
\begin{enumerate}
\item[a.] $f$  is a generalized uniform covering map.
\item[b.] $f$ satisfies Conditions GP1-2.
\item[c.] $\widetilde f\colon GP(X,x_0)\to GP(Y,f(x_0))$
 is a uniform equivalence for some $x_0\in X$.
\end{enumerate}
 \end{theorem}

\proof The equivalence of b) and c) follows from \ref{GUCMInduceUHomeo}.
Suppose $f$ satisfies Conditions GP1-2.
\par {\bf Proof of C1}.
Given an entourage $E$ of $X$ choose an entourage $F_1$ of $X$ so that two generalized paths starting at the same point are $E$-homotopic provided their images are $f(F_1)$-homotopic. Choose an entourage $F$ of $Y$ so that any  $(x,y)\in F$
can be joined by a generalized path that is $f(F_1)$-short. Suppose $(f(x),y)\in F$. Join $f(x)$ and $y$ by a generalized
path $c$ that is $f(F_1)$-short. Now $c$ lifts to a generalized path $\widetilde c$ starting at $x$. Let $y'$ be the
endpoint of $\widetilde c$. Since $c$ is $f(F_1)$-homotopic to the constant path at $f(x)$, $\widetilde c$ is $E$-homotopic to the constant path at $x$. In particular $(x,y')\in E$.

{\bf Proof of C2}.
Given an entourage $E$ of $X$
choose an entourage $G$ of $Y$ so that
$(\widetilde f(c),\widetilde f(d))\in G^\ast$ implies $(c,d)\in E^\ast$ for any
two generalized paths $c$ and $d$ originating from the same point.
Choose an entourage $F\subset H=E\cap f^{-1}(G)$ of $X$ such that any pair $(x,y)\in F$
can be joined by a generalized $H$-short path $c(x,y)$.
Given an $F$-chain $\alpha$ create a generalized uniform path
$p(\alpha)$ by replacing each of its edges $[x_i,x_{i+1}]$ with
$c(x_i,x_{i+1})$. Suppose $f(\alpha)$ is $f(F)$-homotopic
to $f(\beta)$. In that case $(\widetilde f(p(\alpha)),\widetilde f(p(\beta)))\in G^\ast$.
Therefore $(p(\alpha),p(\beta))\in E^\ast$.
As $p(\alpha)_E=[\alpha]$ and $p(\beta)_E=[\beta]$,
$\alpha$ is $E$-homotopic to $\beta$.
\endproof

\begin{corollary}\label{CompositionsOfGenCovers}
Let $X_1$, $X_2$, $X_3$ be uniformly joinable chain connected Hausdorff uniform spaces.
Suppose $f\colon X_1\to X_2$ and $g\colon X_2\to X_3$ generate the uniform structure of their ranges.
If any two of $f$, $g$, $h=g\circ f$ are generalized uniform covering maps, then so is the third.
\end{corollary}

\begin{theorem}\label{piXIsUniversalGenUnifCov}
The projection $\pi_X\colon GP(X,x_0)\to X$ is a universal generalized
uniform covering map in the category of uniformly joinable chain connected
Hausdorff spaces.
\end{theorem}
\proof Suppose $X$ is chain connected and uniformly joinable.
By \ref{GPSquare=GP} and \ref{AltCharOfGenUnifCovForUJ}
$\pi_X$ is a generalized uniform covering map.
If $f\colon X\to Y$ is a generalized uniform covering map,
then for any $x_0\in X$ the induced map $\widetilde f\colon GP(X,x_0)\to GP(Y,f(x_0))$ is a uniform equivalence, so we
put $g=\pi_X\circ \widetilde f^{-1}\colon GP(Y,f(x_0))\to X$. It is clearly a lift of $\pi_Y$.
Since $\widetilde g$ is a uniform equivalence and $GP(Y,f(x_0))$ is joinable, $g$ is a generalized uniform covering map by \ref{AltCharOfGenUnifCovForUJ}.
\endproof

\begin{corollary}\label{GenPathsForUnifPoincare}
If $X$ is a uniform Poincare space, then the natural function
from $\widetilde X$ to $GP(X,x_0)$ is a uniform equivalence.
\end{corollary}
\proof Use \ref{piXIsUniversalGenUnifCov}  to produce
a lift $\alpha\colon GP(X,x_0)\to\widetilde X$ of the projection
$GP(X,x_0)\to X$ that generates the uniform structure of $X$. That lift is the inverse of $\beta\colon \widetilde X\to GP(X,x_0)$ ($\beta$ sends a path in $X$ to the induced generalized uniform path).
Indeed, we can apply \ref{UniquenessOfLifts} to conclude both $\alpha\circ \beta$
and $\beta\circ\alpha$ are identities.
\endproof

\begin{theorem}\label{ExistenceOfLifts}
Suppose $f\colon X\to Y$ is a generalized uniform covering map
and $y_0=f(x_0)$.
If $Z$ is uniformly joinable chain connected, then
the following are equivalent for any $z_0\in Z$ and any uniformly continuous $g\colon Z\to Y$
so that $g(z_0)=y_0$:
\begin{itemize}
\item [a.] There is a uniformly continuous lift $h\colon Z\to X$ of $g$ satisfying $h(z_0)=x_0$,
\item[b.] The image of $\check\pi_1(g)\colon \check\pi_1(Z,z_0)\to \check\pi_1(Y,y_0)$
is contained in the image of $\check\pi_1(f)\colon \check\pi_1(X,x_0)\to \check\pi_1(Y,y_0)$
\end{itemize}
Moreover, if $g$ is a generalized uniform covering map and has a uniformly
continuous lift $h$, then $h$ is a generalized uniform covering map
provided $X$ is joinable.
\end{theorem}
\proof b)$\implies$a).
Given $z\in Z$ pick $c\in GP(Z,z_0)$ from $z_0$ to $z$
and let $d\in GP(X,x_0)$ satisfy $\widetilde f(d)=\widetilde g(c)$.
That $d$ is unique (once $c$ is chosen) and its end-point is our choice for $h(z)$.
If $c^{\prime}$ is another generalized path from $z_0$ to $z$ with the resulting $d^{\prime}\in GP(X,x_0)$,
then $\widetilde g(c\ast (c^{\prime})^{-1})$ is a generalized loop in $Y$ at $y_0$
and we can choose a generalized loop $e\in GP(X,x_0)$ so that $\widetilde f(e)=\widetilde g(c\ast (c^{\prime})^{-1})$.
Now $\widetilde f(e\ast d^{\prime})=\widetilde f(e)\ast \widetilde g(c^{\prime})=\widetilde g(c)=\widetilde f(d)$,
so $e\ast d^{\prime}=d$ and the end-points of $d^{\prime}$ and $d$ are the same.
Hence $h(z)$ is independent on the choice of generalized path $c$.
\par It remains to show $h$ is uniformly continuous and here is where we use
Conditions GP1-2.
Given an entourage $E$ of $X$ pick an entourage $F$ of $Y$ so that
any $F$-short generalized path in $Y$ lifts to an $E$-short generalized path in $X$.
Next choose an entourage $G$ of $Z$ satisfying $g(G)\subset F$.
Finally, choose an entourage $H$ of $Z$ such that any two points
$(z,z^{\prime})\in H$ can be connected by a $G$-short generalized path.
Pick $c\in GP(Z,z_0)$ from $z_0$ to $z$ and then a $G$-short $c^{\prime}$ from $z$
to $z^{\prime}$. The difference between $\widetilde h(c)$ and $\widetilde h(c\ast c^{\prime})$
is $F$-short, so they have lifts to $X$ that differ by an $E$-short path.
The conclusion is that $(h(z),h(z^{\prime}))\in E$ which means $h(H)\subset E$,
i.e. $h$ is uniformly continuous.
\par Assume $X$ is joinable, $g$ is a generalized uniform covering map and has a uniformly
continuous lift $h$. In view of \ref{CompositionOfGenCovers}
it suffices to show $h$ generates the uniform structure of its range.
Since $\widetilde g=\widetilde h\circ \widetilde f$, $\widetilde g$ is a uniform equivalence
and $GP(X,x_0)$ is uniformly joinable.
\endproof

\begin{corollary}\label{GPToXIsACoverInTermsOfpiOne}
Let $X$ be a uniformly joinable and chain connected space.
The projection $\pi_X\colon GP(X,x_0)\to X$ is a uniform covering map
if and only if there is an entourage $E$ of $X$ such that
the natural homomorphism $\check\pi_1(X,x_0)\to \pi_1(R(X,E),x_0)$
is a monomorphism.
\end{corollary}
\proof
Suppose $\pi_X\colon GP(X,x_0)\to X$ is a uniform covering map
and choose
 an entourage $E$ of $X$ such that $E^\ast$ is transverse to $\pi_X$.
 That means
$(c,d)\in E^\ast$ implies $c=d$ if $c$ and $d$ are generalized paths
with the same end-point.
Suppose $c$ and $d$ are generalized loops at $x_0$ so that $c_E=d_E$.
That implies $c^{-1}\ast d$ is $E$-short, hence $c^{-1}\ast d$ is trivial and $c=d$.
\par Suppose the natural homomorphism $\check\pi_1(X,x_0)\to \pi_1(R(X,E),x_0)$
is a monomorphism.
If $c$ and $d$ are two generalized paths from $x_0$ to $x$
such that $b=c^{-1}\ast d$ is $E$-short, then $b_E$ is trivial and $b$ must be trivial.
That means $c=d$ and $\pi_X$ is a uniform covering map.
\endproof

\begin{theorem}\label{GUCMFromMetrizableUJHasCompleteFibers}
Suppose $X$ has a countable base
of entourages and is uniformly joinable.
If $f\colon X\to Y$ is a generalized uniform covering map, then the fibers of $f$ are complete.
\end{theorem}
\proof Choose a base $\{E_n\}_{n=1}^\infty$ of entourages of $X$
satisfying the following conditions:
\begin{enumerate}
\item Any pair $(x,y)\in E_{n+1}$ admits a generalized uniform path
$c_n(x,y)$ from $x$ to $y$ whose $E_n$-term is the edge-path $e(x,y)$.
\item If $\alpha$ and $\beta$ are two $E_{m+1}$-chains
originating at the same point, then they are $E_m$-homotopic
if $f(\alpha)$ is $f(E_{m+1})$-homotopic to $f(\beta)$.
\end{enumerate}

Given a Cauchy sequence in a fiber $f^{-1}(y)$ of $f$ we may choose
its subsequence $\{x_n\}_{n=1}$ such that $(x_k,x_m)\in E_{n+1}$ for $k,m\ge n$.
\par Let $\alpha_1$ be the edge-path $e(x_1,x_2)$.
Given an $E_{n+1}$-chain $\alpha_n$ from $x_1$ to $x_{n+1}$
construct $\alpha_{n+1}$ by replacing each link $e(u,v)$ of $\alpha_n$
by the $E_{n+2}$-term of $c_n(u,v)$ and then concatenating all of it with $e(x_{n+1},x_{n+2})$.
Notice $\{f(\alpha_n)\}_{n=1}$ is a generalized uniform loop at $y$,
so it has a lift $\{\beta_n\}$ from $x_1$ to some $x\in f^{-1}(y)$.
If $x$ is not the limit of $\{x_n\}_{n=1}$, then there is $m\ge 1$
with no $x_i$ belonging to $B(x,E_m)$.
As $f(\alpha_{m+1})$ is $f(E_{m+1})$-homotopic to $f(\beta_{m+1})$,
$\alpha_{m+1}$ is $E_m$-homotopic to $\beta_{m+1}$. In particular,
their end-points are $E_m$-close. Thus $(x,x_{m+1})\in E_m$, a contradiction.
\endproof

\section{Generalized paths relative to spaces}\label{SECTION GenPathsInSpaces}

In this section we expand an idea of Krasinkiewicz and Minc \cite{KraMin}
to define generalized uniform paths of $X$ via an embedding in
a uniform space $T$ with nice local properties.
We require $T$ to be uniformly locally path-connected and the embedding
$X\to T$ satisfies the following analog of uniform semi-local simply connectedness:
Given an entourage $E$ of $T$ there is an entourage $F\subset E$ of $T$
such that any loop in $B(x,F)$ is contractible in $B(X,E)$ for all $x\in X$
(here $B(X,E)$ is $\bigcup\limits_{x\in X}B(x,E)$).

One important case is that of $T$ being uniformly locally simply-connected
as every uniform Hausdorff space $X$ embeds in the Tychonoff cube $I^J$
for some $J$ (that embedding is simply via the set of all
uniformly continuous functions $X\to I$, so that is what one can choose for
index set $J$).
\par Another important case is of $T=X$ and $X$ being a uniform Poincare space.

From now on we assume $X$ is chain-connected.
In this case one can define generalized paths following \cite{KraMin} (only the compact metric case
is discussed there):  $GP_T(X,x_0)$ is the set of {\bf generalized paths} in $T$
from $x_0$ to points of $X$. A {\bf generalized path} is a family $\{[c_E]\}_{E\in \mathcal{E}}$
of homotopy classes of paths $c_E$ in $B(X,E)$ with common end-point $x\in X$
such that for $F\subset E$ the path $c_F$ is homotopic to $c_E$ in $B(X,E)$
rel. end-points.
Given an entourage $F$ of $T$ we define an entourage $F_\ast$ of $GP_T(X,x_0)$
as the set of pairs $(\{[c_E]\}_{E\in \mathcal{E}},\{[d_E]\}_{E\in \mathcal{E}})$
such that $c_F^{-1}\ast d_F$ is homotopic rel. end-points
to a path contained in $B(z,F)$ for some $z\in X$.
\par Notice if $T=X$ and $X$ is a uniform Poincare space,
$GP_X(X,x_0)$ is simply $\widetilde X$.

Our goal is to discuss the connection between $GP(X,x_0)$ and $GP_T(X,x_0)$.

\par Given an entourage $E$ of $T$ let $u(E)\subset E$ be an entourage of $T$
such that any loop in $B(x,u(E)^2)$ is contractible in $B(X,E)$ for all $x\in X$.
Let $v(E)\subset E$ be an entourage of $T$
such that any two points in $B(x,v(E))$ can be connected by a path in $B(x,E)$ for all $x\in T$.
Put $w(E)=v(u(E))$.
\par Given a $w(E)$-chain $c=x_0,\ldots,x_k$ in $X$ from $x_0$ to $x$
choose a path $\alpha_m$ from $x_m$ to $x_{m+1}$ in $B(x_m,u(E))$.
Observe that the homotopy type of $\alpha_m$ (rel. end-points) in $B(X,E)$
does not depend on the choice of $\alpha_m$.
Therefore one has a well-defined path-homotopy class $i(c)$ from $x_0$ to $x$ in $B(X,E)$.

\begin{lemma}\label{IndependenceOfiOnChain}
If $c$ is homotopic to $d$ rel. end-points in $R(X,w(E))$, then $i(c)=i(d)$.
\end{lemma}
\proof It suffices to consider two cases: reduction of a vertex of $x_0,\ldots,x_{k}$
or expansion of $x_0,\ldots,x_{k}$ by a vertex.
\par If a vertex $x_{m+1}$ is dropped from $x_0,\ldots,x_{k}$,
then the concatenation of paths $\alpha_m$ and $\alpha_{m+1}$
is replaced by a path $\beta$ straight from $x_m$ to $x_{m+2}$.
Since $\alpha_m\ast\alpha_{m+1}\ast\beta^{-1}$ is a loop
in $B(x_m,u(E)^2)$, it is null-homotopic in $B(X,E)$
and $\alpha_0\ast\ldots\alpha_{k-1}$ is homotopic rel. end-points
to the concatenation in which $\alpha_m\ast\alpha_{m+1}$ is replaced by $\beta$.
\par The case of expansion of $x_0,\ldots,x_{k}$ by one vertex
is essentially covered by the first case.
\endproof

Given an entourage $F$ of $T$ and given a path $\alpha$ from $x_0$ to $x\in X$
in $B(X,F)$ construct the homotopy class $j(\alpha)$ of a path from $x_0$ to $x$
in $R(X,F^6)$ as follows:
\par For each $t\in [0,1]$ find $x(t)\in X$ so that $(\alpha(t),x(t))\in F$
(obviously, we want $x(0)=x_0$ and $x(1)=x$).
Then find a subdivision $0=t_0\leq t_1\leq\ldots\leq t_k=1$ of the unit interval $I$
such that $\alpha[t_m,t_{m+1}]$ is contained in $B(z,F^2)$ for some $z\in X$.
We need to take $F^2$ since $B(z,F)\subset \mathrm{Int}B(z,F^2)$.
Let $j(\alpha)$ be the homotopy class of the $F^6$-chain
$x(0),\ldots,x(t_k)$ in $R(X,F^6)$.

\begin{lemma}\label{IndependenceOfjOnPath}
$j(\alpha)$ does not depend on the choice of points $0=t_0\leq t_1\leq\ldots\leq t_k=1$
and $j(\alpha)=j(\beta)$ in $R(X,F^{12})$ if $\alpha$ is homotopic to $\beta$ in $B(X,F)$
rel. end-points.
\end{lemma}
\proof To show independence of $j(\alpha)$ of the choice of points $0=t_0\leq t_1\leq\ldots\leq t_k=1$ it suffices to consider the case of expanding $0=t_0\leq t_1\leq\ldots\leq t_k=1$
by adding extra $s$, $t_m\leq s\leq t_{m+1}$. The reason is that any two subdivisions
of the unit interval can be combined by adding one point at the time.
Since $(x(t_m),x(s))\in F^6$ and $(x(t_{m+1}),x(s))\in F^6$,
the chain $x(0),\ldots,x(t_m),x(s),x(t_{m+1}),\ldots,x(t_k)$ is an $F^6$-expansion of
$x(0),\ldots,x(t_k)$ and is homotopic to $x(0),\ldots,x(t_k)$ rel. end-points in $R(X,F^6)$.
\par Suppose $H\colon I\times I\to B(X,F)$ is a homotopy rel. end-points from $\alpha$
to $\beta$. There is an equally spaced subdivision $0=t_0\leq t_1\leq\ldots\leq t_k=1$ of the unit interval $I$
so that $H([t_m,t_{m+1}]\times [t_n,t_{n+1}])\subset B(z,F^2)$ for some $z\in X$.
To conclude $j(\alpha)=j(\beta)$ in $R(X,F^{12})$ it suffices to apply the following:
\par {\bf Observation}. If $E$ is an entourage of $X$ and $x_0,\ldots,x_k$,
$y_0,\ldots,y_k$ are two $E$-chains joining $x_0$ to $x$, then
they are homotopic in $R(X,E^2)$ rel. end-points if $(x_n,y_n)\in E$ for all $n\leq k$.
\par {\bf Proof of Observation}.
Create an $E^2$-chain $x_0,y_0,\ldots, x_k,y_k$ and notice it can be reduced to both
$x_0,\ldots,x_k$ and $y_0,\ldots,y_k$ in $R(X,E^2)$.\endproof

\begin{lemma}
\label{ijAndjiAreInclusion} Let $E$ be an entourage of $T$ and $F$ be an
entourage with $F^{12}\subset w(E)$. If $\alpha $ is a path in $B(X,F)$ then
$i(j(\alpha ))$ is homotopic to $\alpha $ in $B(X,E)$. Similarly, let $E$ be
an entourage of $T$ and $F$ be an entourage with $F^{12}\subset E$. If $
\gamma $ is a path in $R(X,w(F))$, then $j(i(\gamma ))$ is homotopic to $
\gamma $ in $R(X,E)$.
\end{lemma}

\proof Say $j(\alpha )$ is the homotopy class of $
x_{0}=x(t_{0}),\ldots ,x(t_{k})$. For each $i<k$, $\alpha (t_{i}),\alpha
(t_{i+1}),x(t_{i+1})\in B(x(t_{i}),w(E))$ so there are paths from $x(t_{i})$
to $\alpha (t_{i})$ and from $x(t_{i+1})$ to $\alpha (t_{i+1})$ that are
contained in $B(x(t_{i}),u(E))$. Therefore $i(j(\alpha ))$ is homotopic to $
\alpha $ in $B(X,E)$. Now suppose $\gamma $ is represented by the $w(E)$-chain $x_{0},\ldots x_{k}$. Notice $j(i(\gamma ))$ is the same chain in $R(X,E)$
since for each $i$, $\alpha _{i}$ (from the definition of $j)$ is contained
in $B(x_{i},u(F))\subset B(x_{i},F^{2})$.\endproof

Now we are in a position to define $i\colon GP(X,x_0)\to GP_T(X,x_0)$
and  $j\colon GP_T(X,x_0)\to GP(X,x_0)$.

Given $c=\{[c_F]\}_{F\in\mathcal{E}}\in GP(X,x_0)$ (from $x_0$ to $x$) assume each $c_F$
is realized by an $F$-chain $x_0,\ldots,x_{k(F)}$ in $R(X,F)$.
Given an entourage $E$ of $T$ use \ref{IndependenceOfiOnChain}
to notice $i(c_F)$ is independent of the choice of $F\subset w(E)$.
By putting $i(c)_E=i(c_F)$ we get a well-defined element of $GP_T(X,x_0)$.
Similarly, given an element $\alpha=\{[\alpha_F]\}_{F\in\mathcal{E}}\in GP_T(X,x_0)$
use \ref{IndependenceOfjOnPath} to notice that for any entourage $E$ of $T$
the element $j(\alpha_F)$ does not depend on $F$ provided $F^{12}\subset E$.
Thus, putting $j(\alpha)=\{j(\alpha_F)\}_{E\in\mathcal{E}}$ we get a well-defined element
of $GP(X,x_0)$.

\begin{theorem}
\label{iFromGPToGP_TIsUnifCont} $i\colon GP(X,x_{0})\rightarrow
GP_{T}(X,x_{0})$ and $j\colon GP_{T}(X,x_{0})\rightarrow GP(X,x_{0})$ are
uniformly continuous and inverse to each other.
\end{theorem}

\proof For uniform continuity of $i$ let us show that $i(w(E)^{\ast
})\subset E_{\ast }$. Let $(c,d)\in w(E)^{\ast }$. Then $c_{w(E)}^{-1}\ast
d_{w(E)}$ is homotopic to $e(x,y)$ in $R(X,w(E))$ where $x$ and $y$ are the
endpoints of $c$ and $d$ respectively. Then we have $i(c)_{E}^{-1}\ast
i(d)_{E}=i(c_{w(E)})^{-1}\ast i(d_{w(E)})=i(c_{w(E)}^{-1}\ast
d_{w(E)})=i(e(x,y))$ so $(i(c),i(d))\in E_{\ast }$. Similarly, for
continuity of $j$, we show $j(F_{\ast })\subset E^{\ast }$ for $
F^{12}\subset E$. Let $(\alpha ,\beta) \in F_{\ast }$ and $x$ and $y$ be the
endpoints of $\alpha $ and $\beta $ respectively. Then $\alpha _{F}^{-1}\ast
\beta _{F}$ is homotopic in $B(X,F)$ to some path $\gamma $ from $x$ to $y$
that is contained in $B(z,F)$ for some $z\in X$. Then we have $j(\alpha
)_{E}^{-1}\ast j(\beta )_{E}=j(\alpha _{F})^{-1}\ast j(\beta _{F})=j(\alpha
_{F}^{-1}\ast \beta _{F})=j(\gamma )$ in $R(X,E)$. Notice that in $R(X,E)$, $
j(\gamma )$ is homotopic to $e(x,y)$.

Let $\alpha \in GP_{T}(X,x_{0})$ and consider $i(j(\alpha ))$. We have $
i(j(\alpha ))_{E}=i(j(\alpha )_{w(E)})=i(j(\alpha _{F}))$ where $
F^{12}\subset w(E)$. By \ref{ijAndjiAreInclusion}, $i(j(\alpha _{F}))$ is
homotopic to $\alpha _{E}$ in $B(X,E)$. Now let $c\in GP(X,x_{0})$ and
consider $j(i(c))$. We have $j(i(c))_{E}=j(i(c)_{F})=j(i(c_{w(F)}))$ for $
F^{12}\subset E$. Again, by \ref{ijAndjiAreInclusion}, $j(i(c_{w(F)}))$ is
homotopic in $R(X,E)$ to $c_{E}$.\endproof

 \begin{corollary}\label{FirstShapeGroupAndDeckGroup}
 If $X$ is a metric continuum, then $\check\pi_1(X,x_0)$
 is isomorphic to the first shape group of $(X,x_0)$.
\end{corollary}
\proof Embed $X$ in the Hilbert cube $Q$.
As in the proof of \ref{iFromGPToGP_TIsUnifCont}, $\check\pi_1(X,x_0)$
is isomorphic (also in the category of topological groups)
to the group of generalized loops of $X$ in $Q$ at $x_0$.
That is the same as
the inverse limit of $\{\pi_1(U_n)\}$ (each with the discrete topology),
where $U_n$ is the $\frac{1}{n}$-ball of $X$ in $Q$,
and that is exactly the first shape group of $(X,x_0)$ (see \cite{DydSeg} or
\cite{MarSeg}).
\endproof

\begin{definition}[cf. {\cite[p.88]{DydSeg}}]
A pointed continuum $(X,x_0)$ is called {\bf pointed $1$-movable} if $\mathrm{pro-}\pi_1(X,x_0)$
satisfies the Mittag-Leffler condition.
\end{definition}

\begin{corollary}\label{CharOfJoinableContinua}
For a metric continuum $X$ the following conditions are equivalent:
\begin{itemize}
\item[a.] $X$ is joinable,
\item[b.] $X$ is pointed $1$-movable.
\item[c.] $X$ is uniformly joinable.
\end{itemize}
\end{corollary}
\proof
a)$\implies$b).
Embed $X$ in the Hilbert cube $Q$.
As in the proof of \ref{iFromGPToGP_TIsUnifCont} joinability of $X$
is equivalent to the property that every two points $x,y\in X$
there is a sequence of paths $a_n$ joining $x$ to $y$ in the $(1/n)$-neighborhood
$U_n$ of $X$ such that $a_{n+1}$ is homotopic to $a_n$ rel. end-points in
$U_n$. That coincides with the original definition of
joinability of $X$ given by Krasinkiewicz and Minc \cite{KraMin}.
The main result of \cite{KraMin} states that joinable continua
have the fundamental pro-group satisfying the Mittag-Leffler condition.
\par b)$\implies$c).
By~\cite[Theorem 6.1.7]{DydSeg}, we have ${\lim\limits_{\leftarrow}}^1 (\mathrm{pro-}\pi_1(X,x_0))=0$.
By~\ref{Lim1AndJoinability} (or \cite{KraMin}), $X$ is joinable.
Applying~\ref{countable ML implies strong ML}
and~\ref{StrongMLPlusJoinImpliesUJ}
gives $X$ being uniformly joinable.
c)$\implies$a).
Follows from~\ref{ChainConnectedAndUJAreJoinable} since any metric continuum is chain connected.
\endproof

Since all subcontinua of surfaces are pointed $1$-movable (see \cite[Theorem 7.1.7]{DydSeg}),
one has the following:

\begin{corollary}\label{ExamplesOfPointedOneMovable}
All subcontinua of surfaces are uniformly joinable (that includes the suspension of the
Cantor set and the Hawaiian Earring). The dyadic solenoid is not joinable.
\end{corollary}

In connection to \ref{ExamplesOfPointedOneMovable}
let us point out the boundary of any word-hyperbolic group is compact and
metrizable~\cite{Gromov} and the boundary of any one-ended word-hyperbolic group is locally
connected~\cite{Bowditch} (hence pointed $1$-movable).
Also, pointed $1$-movability is related to
semi-stability at infinity of groups (see \cite{MihTsc} and \cite{Geoghegan}).

\begin{corollary}
\label{GPToXIsACoverInTermsOfpiOneForContinua} If $X$ is a uniformly
joinable metric continuum, then the following conditions are equivalent:

\begin{itemize}
\item[a.] The projection $\pi_X\colon GP(X,x_0)\to X$ is a
uniform covering map,

\item[b.] $\check{\pi}_{1}(X,x_{0})$ is countable,

\item[c.] $\check{\pi}_{1}(X,x_{0})$ is finitely generated.
\end{itemize}
\end{corollary}

\proof Embed $X$ in the Hilbert cube $Q$. We show that $\pi
_{X}:GP(X,x_{0})\rightarrow X$ is a uniform covering map if and
only if there is a closed neighborhood $N$ of $X$ in $Q$ with $\check{\pi}
_{1}(X,x_{0})\rightarrow \pi _{1}(N,x_{0})$ a monomorphism and $N$ the
homotopy type of a compact polyhedron. That condition is known to be
equivalent to b) and c) (see \cite{Dyd4} or \cite[Corollary 8 on p.177]{MarSeg}). According to~\ref{GPToXIsACoverInTermsOfpiOne} $\pi _{X}:GP(X,x_{0})\rightarrow X$ is a
uniform covering map if and only if there is an entourage $E$ of $X$
so that $\check{\pi}_{1}(X,x_{0})\rightarrow \pi _{1}(R(X,E),x_{0})$ is a
monomorphism. Suppose such an $E$ exists and let $F$ be an entourage of $X$
with $F^{12}\subset E$. Then $\check{\pi}_{1}(X,x_{0})\rightarrow \pi
_{1}(R(X,w(F)),x_{0})\rightarrow \pi _{1}(B(X,F),x_{0})\rightarrow \pi
_{1}(R(X,E),x_{0})$ is a monomorphism (see \ref{ijAndjiAreInclusion}) so $
\check{\pi}_{1}(X,x_{0})\rightarrow \pi _{1}(B(X,F),x_{0})$ is as well. Note
that there is a closed neighborhood $N$ of $X$ in $Q$ with $N\subset B(X,F)$
and $N$ the homotopy type of a compact polyhedron. Similarly, if such an $N$
exists, Find $\varepsilon >0$ so that $B(X,E_{\varepsilon })\subset N$ where
$E_{\varepsilon }=\{(x,y)\in Q\times Q:d(x,y)<\varepsilon \}$. Take an
entourage $F$ such that $F^{12}\subset w(E_{\varepsilon })$. Then $\check{\pi
}_{1}(X,x_{0})\rightarrow \pi _{1}(B(X,F),x_{0})\rightarrow \pi
_{1}(R(X,w(E_{\varepsilon })),x_{0})\rightarrow \pi _{1}(B(X,E_{\varepsilon
}))$ is a monomorphism so $\check{\pi}_{1}(X,x_{0})\rightarrow \pi
_{1}(R(X,w(E_{\varepsilon })),x_{0})$ is as well.\endproof

\begin{theorem}\label{UniversalCoverAndGUCM}
Suppose $X$ is a path-connected uniform space.
If the projection $\pi\colon \widetilde X\to X$ is a generalized uniform covering map, then
$X$ is a uniform Poincare space and $\pi$ is a uniform covering map.
\end{theorem}
\proof $X$ is uniformly locally path-connected by \ref{UnifLocalPathConnectCharacterization}.
It suffices to show $X$ is uniformly semi-locally simply connected.
Suppose for each entourage $E$ of $X$ there is a point $x_E$
and a loop $\alpha_E$ at $x_E$ in $B(x,E)$ that is non-trivial in $X$.
Pick a path $\gamma_E$ from $x_0$ to $x_E$.
By picking points on the loop $\beta_E=\gamma_E\ast \alpha_E\ast \gamma_E^{-1}$
that belong to the image of $\gamma_E$ only
one can define an $E^\ast$-chain in $\widetilde X$ starting from the trivial loop at $x_0$
and ending at $\beta_E$.
The same chain works for $\omega_E=\gamma_E\ast \gamma_E^{-1}$ but this time we
do not go around $\alpha_E$. Thus we have two $E^\ast$-chains in $\widetilde X$
with the same projection in $X$, so they should be $X\times X$-homotopic
for some $E$, a contradiction.
\endproof

\section{Comparison to Berestovskii-Plaut uniform covers}\label{SECTION-Comparison}

Berestovskii and Plaut used an analogue of the the Schreier
construction for topological groups \cite{BP1} to create an inverse
limit construction \cite{BP3} for a uniform space $X$. We recall
their construction (which we denote by ${\widetilde X}_{BP}$
as $\widetilde X$ is used by us for classical universal cover) below, and compare their inverse limit space
${\widetilde X}_{BP}$ to $GP(X,x_0)$.

Let $X$ be an uniform space with a fixed base point $x_0$. For any
entourage $E$ an $E$-chain starting at $x_0$ and ending at $x \in X$
is a finite sequence of points $\{x_0,\ldots,x_n=x\}$ such that
$(x_i,x_{i+1}) \in E$ for all $ 0 \le i \le n-1$. An $E$-extension
of a $E$-chain $\{x_0,\ldots,x_n=x\}$ is a $E$-chain
$\{x_0,\ldots,x_i,y,x_{i+1},\ldots,x_n=x\}$, with $0 \le i<n$. An
$E$-homotopy is a finite sequence of $E$-extensions (or their
obvious analogues $E$-contractions). $X_E$ is the set of all
$E$-homotopy classes $[c]_E$ of $E$-chains $c$. For any entourage $F
\subset E$ define $\hat{F}$ as follows: $([c]_E,[d]_E) \in \hat{F}$
if $([c]_E,[d]_E)=([x_0,\cdots,x_n,y]_E,[x_0,\cdots,x_n,z]_E)$ with
$ (y,z) \in F$. The collection of all such $\hat{F}$ is a base for
the uniformity on $X_E$. If $F \subset E$ is an entourage, there is
a natural map ${\phi}_{EF}:X_F \to X_E$
which sends $[c]_F$ to $[c]_E$ and generates the uniform structure of $X_E$. 
With hindsight one may say the structure on $X_E$ mimicks the basic topology on $\widetilde X$.

The inverse limit ${\widetilde
X}_{BP}$ of $\{X_E\}_{E\in \mathcal{E}}$ is given the inverse limit uniformity.
Thus ${\widetilde X}_{BP}$ is equivalent to our space $GP(X,x_0)$.
The advantage of our description is a closer connection to the classical
universal cover $\widetilde X$ and generalized paths of Krasinkiewicz-Minc.
\par For the same reason the deck group $\delta_1(X)$ of \cite{BP3}
is isomorphic to our fundamental uniform group $\check\pi_1(X,x_0)$.
Again, the advantage of our approach is the connection between
$\check\pi_1(X,x_0)$ and the fundamental shape group in case
of metrizable compact spaces $X$.
\par The basic class of uniform spaces for which the approach in
\cite{BP3} works is the class of coverable spaces.
A uniform space $X$ is {\bf coverable} if there is a
uniformity base of entourages $E$ (including $X \times X$) such that
the projections ${\widetilde X}_{BP} \to X_E$ are surjective.
In our language that means for every path $\alpha$ in $R(X,E)$
there is a generalized path $c=\{c_F\}_{F\in\mathcal{E}}$ such that
$c_E$ is homotopic rel.end-points to $\alpha$. Thus every coverable
space is uniformly joinable and our theory of generalized uniform covering maps
induces most basic results of \cite{BP3}.
A natural question arises:
\begin{problem}\label{CoverableVsUniformlyJoinable}
Is every uniformly joinable chain-connected space coverable?
\end{problem}

The relevance of \ref{CoverableVsUniformlyJoinable} is that it would
imply a positive answer to Problem
106 of \cite{BP3} for chain-connected spaces (that problem asks if $X$ is coverable provided
${\widetilde X}_{BP} \to X$ is a
uniform equivalence). Indeed, \ref{BiUniformAndJoinability} implies
$X$ is uniformly joinable if ${\widetilde X}_{BP} \to X$ is a
uniform equivalence.

\par There are two obvious strategies to solve \ref{CoverableVsUniformlyJoinable} positively:
\begin{itemize}
\item[a.] Given an entourage $E$ of $X$ choose an entourage $F\subset E$
with the property that any pair $(x,y)\in F$ can be
connected by a generalized path $c(x,y)$ so that its $E$-term is the edge $e(x,y)$.
Try to show ${\widetilde X}_{BP} \to X_F$ is surjective.
\item[b.] Given an entourage $E$ of $X$ define $G(E)$ as all pairs $(x,y)\in E$
with the property that there are generalized paths $c$ from $x_0$ to $x$ and $d$
from $x_0$ to $y$
such that $(c^{-1}\ast d)_E$ is homotopic in $R(X,E)$ to the edge $e(x,y)$
(as $X$ is uniformly joinable $G(E)$ contains $F$ above and is an entourage).
Try to show ${\widetilde X}_{BP} \to X_{G(E)}$ is surjective.
\end{itemize}

Notice Strategy b) is a natural reaction once one realizes Strategy a) fails.

Let us show two examples negating the above strategies.

\begin{example}\label{BasicExampleForBadStrategy}
Consider a regular hexagon with one edge $ab$ of size $1$ removed.
Let $E$ be pairs of distance at most $3$ and $F$ are pairs of distance at most $1$.
\end{example}
\proof
To check that any $F$-short pair can be connected by the right path (notice there is only
one path for every pair anyhow) it suffices
to prove it for $(a,b)$. Let $\alpha$ be the genuine path in $X$ from $a$ to $b$.
We can eliminate first all non-vertex points,
then all vertices and $\alpha_E$ is homotopic in $R(X,E)$ to $e(a,b)$.
Here is the problem: consider chain $x_0=a, x_1=b$ in $F$
and suppose there is a generalized path $\alpha$ whose $F$-term
is homotopic to $\{x_0,x_1\}$. There is no way a $1$-chain from $a$ to $b$
to be $1$-homotopic to $\{x_0,x_1\}$ (consider the last point removed prior to arriving
at pair $\{x_0,x_1\}$) and such generalized path would produce a chain of that kind.
\endproof

\begin{example}
Consider a regular hexagon with one edge $ab$ of size $1$ removed. Add the center $c$ of the hexagon
plus a vertical regular hexagon with bottom $ac$ that we remove. The resulting $X$
and $E=\{(x,y) | dist(x,y)\leq 1=dist(a,b)\}$ have the property
that $(a,b)\in G(E)$ but they cannot be joined by a generalized path
in $X$ whose $G(E)$-term is the edge as $(p,c)\notin G(E)$
for any point $p$ belonging to the first hexagon.
\end{example}

Example \ref{BasicExampleForBadStrategy} says there is an error
in \cite{Platfus}. Indeed, in the proof of Proposition 5 one considers
the entourage $F^\ast$  in $X_E$ consisting of pairs of homotopy classes of paths
$(a,b)$ such that their end-points $x$ and $y$ satisfy $(x,y)\in F$,
$a^{-1}\ast b$ is homotopic rel.end-points to the edge $e(x,y)$,
and there are generalized paths $c$ and $d$ so that $c_E=a$
and $d_E=b$. The entourage $G$ in $\widetilde X=\widetilde X_{BP}$ is defined
as pairs $(c,d)$ so that $(c_E,d_E)\in F^\ast$ and Proposition 5
claims the projection $\widetilde X_G\to \widetilde X$ is a homeomorphism for all such $G$.
Once that holds the proof of Lemma 6 in \cite{Platfus} gives that $\widetilde X\to X_{\pi(G)}$
is surjective provided all such defined entourages $G$ form a base of entourages of $\widetilde X$
which is so if $X$ is uniformly joinable.
However, $\pi(G)$ ($\pi$ being the projection from $\widetilde X$ to $X$)
is exactly $F$ and Example \ref{BasicExampleForBadStrategy} shows
the projection $\widetilde X\to X_{\pi(G)}$ may not be surjective.

The best way to explain to a topologist the philosophical difference between Berestovskii-Plaut
notion of coverability and our notion of uniform joinability
is to point out the latter is a $UV$-type condition and the former one is the same condition
replaced by existence of a base where $V$ can be chosen equal to $U$.
In Siebenmann's thesis he starts from $UV$-type conditions and produces
an end of a manifold. Such an end can be intuitively explained by requiring $V=U$
for some base of neighborhoods $U$ of infinity and some $UV$-type condition.
That means an answer to \ref{CoverableVsUniformlyJoinable} could be positive
but a topologist would be sceptical without adding extra conditions
on the space $X$.

From algebraic point of view uniform joinability corresponds to the Mittag-Leffler condition
and, for inverse sequences of groups, Mittag-Leffler condition is indeed equivalent
to existence of an inverse sequence of epimorphisms. That analogy may lead
to a larger dose of optimism in a positive answer to \ref{CoverableVsUniformlyJoinable}.
However, one may point out that Theorem 7 of \cite{BP2}
characterizes coverability of a locally compact topological group $G$ as being equivalent
to $G$ being connected and locally arcwise connected.
Thus, \ref{CoverableVsUniformlyJoinable} has a positive answer for locally compact
topological groups which may be analogous to the Mittag-Leffler condition
for inverse sequences of groups.

Summing up: uniform joinability is of a shape-theoretical
nature and coverability is more of a geometrical nature.

In \cite{BP3} (on p.1751, the paragraph below Theorem 3) the authors
mention they do not know whether the composition of pro-discrete covers between coverable spaces (or uniform spaces in
general) is a pro-discrete cover but in the case of topological groups it is so.
Let us point out an example resolving that question in the negative even for discrete covers.

\begin{example}\label{CompositionExample}
Consider the case of subgroups $G_1\subset G_2$ of a group $G_3$
such that $G_2$ is normal in $G_3$, $G_1$ is normal in $G_2$ but $G_1$
is not normal in $G_3$.
In case of $G_1=Z$, the group of integers, $G_2=Z\times Z$,
and $G_3$ as the HNN-extension of $G_2$ that switches the $Z$-factors,
the corresponding space is a Seifert 3-manifold which is a locally trivial
fibration over a circle with the fiber homeomorphic to a torus such that
the monodromy (along the base circle) is a homeomorphism of the fiber that
switches the meridian and the parallel of the fiber.

Choose a simplicial complex $K$ with $\pi_1(K)=G_3$, create a covering
$p\colon L\to K$ with $\pi_1(L)=G_2$ and $\pi_1(p)$ realizing inclusion $G_2\to G_3$.
Similarly, create a covering $q\colon M\to L$ with $\pi_1(M)=G_1$
and $\pi_1(q)$ realizing inclusion $G_1\to G_2$.
Notice we can give $L$ and $M$ structures of simplicial complexes (with resulting uniform structures generated by simplicial metrics)
so that both $p$ and $q$ are simplicial maps. Obviously, $p\circ q$
is a simplicial covering map. However, it cannot be realized as a result
of an equi-continuous action of any group $G$ on $M$.
Indeed, as $G_1$ is not normal in $G_3$ there is a loop $\alpha$ in $K$
with two lifts $\beta$ and $\gamma$ (originating at different points $x$ and $y$ of $K$, obviously)
such that $\beta$ is a loop and $\gamma$ is not a loop.
Choose $g\in G$ so that $g\cdot x=y$. Notice $g\cdot \alpha$ and $\gamma$
are two different lifts of the same loop, a contradiction.
\end{example}

\end{document}